\documentclass[11pt,amsfonts]{article}
\usepackage{graphicx}
\usepackage{latexsym}
\usepackage{amssymb}
\usepackage{amsmath}
\usepackage{layout}
\newtheorem{prop}{Proposition}
\newtheorem{lemma}{Lemma}

\newtheorem{theorem}{Theorem}
\newtheorem{remark}{Remark}

\def\real{{\mathord{{\rm I\kern-2.8pt R}}}}        % Fake blackboard bold R.
\def\inte{{\mathord{{\rm I\kern-2.8pt N}}}}

\def\sZZ{{\rm Z\kern-2.8ptem{}Z}}

\def\z{{\mathchoice
  {\sZZ}
  {\sZZ}
  {\rm Z\kern-0.30em{}Z}
  {\rm Z\kern-0.25em{}Z} }}
\def\sQQ{{\kern 0.27em \vrule height1.45ex width0.03em depth0em
          \kern-0.30em \rm Q}}
\def\qu{{\mathchoice
    {\sQQ}
    {\sQQ}
  {\kern 0.225em \vrule height1.05ex width0.025em depth0em \kern-0.25em \rm Q}
  {\kern 0.180em \vrule height0.78ex width0.020em depth0em \kern-0.20em \rm Q}
        }}
\def\sCC{{\kern 0.27em \vrule height1.45ex width0.03em depth0em
          \kern-0.30em \rm C}}
\def\complex{{\mathchoice
    {\sCC}
    {\sCC}
  {\kern 0.225em \vrule height1.05ex width0.025em depth0em \kern-0.25em \rm C}
  {\kern 0.180em \vrule height0.78ex width0.020em depth0em \kern-0.20em \rm C}
        }}

\newenvironment{dem}{{\bf Proof:}}{}

\newcommand{\ba}{\begin{array}}
\newcommand{\ea}{\end{array}}
\newcommand{\be}{\begin{equation}}
\newcommand{\ee}{\end{equation}}
\newcommand{\bea}{\begin{eqnarray}}
\newcommand{\eea}{\end{eqnarray}}
\newcommand{\beaa}{\begin{eqnarray*}}
\newcommand{\eeaa}{\end{eqnarray*}}

\newcommand{\eps}{\varepsilon}

\def\z{\zeta}

\font\tenmath=msbm10 \font\sevenmath=msbm7 \font\fivemath=msbm5
\newfam\mathfam \textfont\mathfam=\tenmath
\scriptfont\mathfam=\sevenmath \scriptscriptfont\mathfam=\fivemath

\def \={{\buildrel {\rm (law)} \over =}}

\def\qed{ \hfill \vrule width.25cm height.25cm depth0cm\smallskip}

\newcommand{\basa}{\begin{assumption}}
\newcommand{\easa}{\end{assumption}}

\newcommand{\bas}{\begin{assum}}
\newcommand{\eas}{\end{assum}}

\newcommand{\ignore}[1]{}
\begin{document}

\renewcommand{\thefootnote}{\fnsymbol{footnote}}

\renewcommand{\thefootnote}{\fnsymbol{footnote}}

\date{  }
\title{Asymptotic theory for fractional regression models via Malliavin calculus}
\author{ Solesne Bourguin $^{1}\qquad $%
Ciprian A. Tudor $^{2,}$ \footnote{Associate member of the team Samm, Universit\'e de Panth\'eon-Sorbonne Paris 1 }\vspace*{0.1in} \\
%\thanks{\emph{Corresponding author}; research
%partially supported by NSF grant DMS 0606615.}
$^{1}$SAMM, Universit\'e de Paris 1 Panth\'eon-Sorbonne,\\
90, rue de Tolbiac, 75634, Paris, France. \\
solesne.bourguin@univ-paris1.fr \vspace*{0.1in} \\
 $^{2}$ Laboratoire Paul Painlev\'e, Universit\'e de Lille 1\\
 F-59655 Villeneuve d'Ascq, France.\\
 \quad tudor@math.univ-lille1.fr\vspace*{0.1in}}
\maketitle

\begin{abstract}
\noindent We study the asymptotic behavior as $n\to \infty$ of the sequence
$$S_{n}=\sum_{i=0}^{n-1} K(n^{\alpha} B^{H_{1}}_{i}) \left( B^{H_{2}}_{i+1}-B^{H_{2}}_{i}\right)$$
where $B^{H_{1}}$ and $B^{H_{2}}$ are two independent fractional Brownian motions, $K$ is a kernel function and the bandwidth parameter $\alpha$ satisfies certain hypotheses in terms of $H_{1}$ and $H_{2}$. Its limiting distribution is a mixed normal law involving the local time of the fractional Brownian motion $B^{H_{1}}$.
We use the techniques of the Malliavin calculus with respect to the fractional Brownian motion.

\end{abstract}

\vskip0.3cm

{\bf  2010 AMS Classification Numbers:}  60F05, 60H05, 91G70.

 \vskip0.3cm


\section{Introduction}
The motivation of our work comes from the econometric theory. Consider a regression model of the form
\begin{equation*}
y_{i}= f(x_{i}) +u_{i}, \hskip0.5cm i\geq 0
\end{equation*}
where $(u_{i})_{i\geq 0}$ is the "error" and $(x_{i})_{i\geq 0}$ is the regressor. The purpose is to estimate the function $f$ based on the observation of the random variables $y_{i}$, $i\geq 0$. The conventional kernel estimate of $f(x)$ is
\begin{equation*}
\hat{f}(x)= \frac{ \sum_{i=0}^{n} K_{h} (x_{i}-x) y_{i} }{\sum_{i=0}^{n} K_{h}(x_{i}-x)}
\end{equation*}
where $K$ is a nonnegative real kernel function satisfying  $\int_{\mathbb{R}}K^{2}(y)dy =1$ and $\int_{\mathbb{R}}yK(y)dy =0$
and $K_{h}(s)=\frac{1}{h}K(\frac{s}{n})$. The bandwidth parameter $h\equiv h_{n}$ satisfies $h_{n}\to 0$ as $n\to \infty$.
The asymptotic behavior of the estimator $\hat{f}$ is usually related to the behavior of the sequence
\begin{equation*}
V_{n}= \sum_{i=1}^{n} K_{h} (x_{i}-x) u_{i}.
\end{equation*}
The  limit in distribution as $n\to \infty$ of the sequence $S_{n}$ has been widely studied in the literature in various situations. We refer, among others, to \cite{KT} and \cite{KMT} for the case where $x_{t}$ is a recurrent Markov chain, to \cite{WP1} for the case where $x_{t}$ is a partial sum of a general linear process, and \cite{WP2} for a more general situation. See also \cite{PP} or \cite{Ph}. An important assumption in the main part of the above references is the fact that $u_{i}$ is a martingale difference sequence. In our work we will consider the following situation: we assume that the regressor $x_{i}= B^{H_{1}}_{i}$ is a fractional Brownian motion (fBm) with Hurst parameter $H_{1}\in (0,1)$ and the error is $u_{i}=B^{H_{2}}_{i+1}-B^{H_{2}}_{i}$  where $B^{H_{2}}$ is a fBm with $H_{2}\in (0,1)$ and  it is independent from $B^{H_{1}}$. In this case, our error process has no semimartingale property. We will also set $h_{n}=n^{-\alpha}$ with $\alpha >0$. A supplementary assumption on $\alpha$ will be imposed later in terms of the Hurst parameters $H_{1}$ and $H_{2}$. The sequence $V_{n}$ can be now written as
\begin{equation}\label{snx}
S_{n}(x)= \sum_{i=0} ^{n} K(n^{\alpha} (B^{H_{1}}_{i} - x))\left( B^{H_{2}}_{i+1}-B^{H_{2}}_{i}\right).
\end{equation}
Our purpose is to give an approach based on stochastic calculus for this asymptotic theory. Recently, the stochastic integration with respect to the fractional Brownian motion has been widely studied. Various types of stochastic integrals, based on Malliavin calculus, Wick products or rough path theory have been introduced and change of variables formulas have been derived. We will use all these different techniques in our work. The general idea is as follows. Suppose that $x=0$. We will first observe that the asymptotic behavior of the sequence $S_{n}$ will be given by the sum
\begin{equation}\label{an}
a_{n}=\sum_{i,j=0} ^{n}  K(n^{\alpha} B^{H_{1}}_{i})K(n^{\alpha} B^{H_{1}}_{j})\mathbf{E}\left( (B^{H_{2}}_{i+1}-B^{H_{2}}_{i})(B^{H_{2}}_{j+1}-B^{H_{2}}_{j})\right).
\end{equation}
This is easy to understand since the conditional distribution of $S_{n}$ given $B^{H_{1}}$ is given by
$$\left( a_{n}\right) ^{\frac{1}{2}} Z$$
where $Z$ is a standard normal random variable. The double sum  $a_{n}$ can be decomposed into two parts: a ``diagonal'' part given by $\sum_{i=1}^{n}K^{2}(n^{\alpha}B^{H_{1}}_{i})$ and a ``non-diagonal'' part given by the terms with $i\not=j$. We will restrict ourselves to the situation where the diagonal part is dominant (in a sense that will be defined later) with respect to the non-diagonal part. This will imply a certain assumption on the bandwidth parameter $\alpha$ in terms of $H_{1}$ and $H_{2}$. We will therefore need to study the asymptotic behavior of
\begin{equation}\label{brac}
\langle S\rangle _{n}:= \sum_{i=1}^{n}K^{2}(n^{\alpha}B^{H_{1}}_{i}).
\end{equation}
(In the case $H_{2}=\frac{1}{2}$  this is actually the bracket of  $S_{n}$  which is a martingale; this motivates our choice of notation.) We will assume that the kernel $K$ is the standard Gaussian kernel
\begin{eqnarray}
K(x) = \frac{1}{\sqrt{2\pi}}e^{-\frac{x^{2}}{2}}. \nonumber
\end{eqnarray}
This choice is motivated by the fact that $K^{2}(n^{\alpha}B^{H_{1}}_{i})$ can be decomposed into an orthogonal sum of multiple Wiener-It\^o integrals (see \cite{NV}, \cite{CNT}, \cite{Edd}) and the Malliavin calculus can be used to treat the convergence of (\ref{brac}). Its limit in distribution will be after normalization the local time of the fractional Brownian motion denoted $cL^{H_{1}}(1,0)$, where $c$ is positive constant. Consequently, we will find that the (renormalized) sequence $S_{n}$ converges in law to a mixed normal random variable $cW_{L^{H_{1}}(1,0)}$ where $W$ is a Brownian motion independent from $B^{H_{1}}$ and $c$ is a positive constant. The result is in concordance with the papers \cite{WP1}, \cite{WP2}.
\\\\
But we also prove a stronger result: we show that the vector $(S_{n}, (G_{t})_{t\geq 0} )$ converges in the sense of finite dimensional distributions to the vector $(cW_{L^{H_{1}}(1,0)}, (G_{t})_{t\geq 0})$, where $c$ is a positive constant, for any  stochastic process $(G_{t})_{t\geq 0}$ independent from $B^{H_{1}}$ and adapted to the filtration generated by $B^{H_{2}}$ which satisfies some regularity properties in terms of the Malliavin calculus. We will say that $S_{n}$ converges stably to its limit. To prove this stable convergence we will express $S_{n}$ as a stochastic integral with respect to $B^{H_{2}}$ and we will use the techniques of the Malliavin calculus. We will limit ourselves in this last section to the case $H_{2}>\frac{1}{2}$.
\\\\
We also mention that, although the error process $B^{H_{2}}$ does not appear in the limit of (\ref{snx}), it governs the  behavior of this sequence. Indeed, the parameter $H_{2}$ is involved in the renormalization of (\ref{snx}) and the stochastic calculus with respect to $B^{H_{2}}$ is crucial in the proof of our main results.
\\\\
We have organized our paper as follows: Section 2 contains the notations, definitions and results from the stochastic calculus that will be needed throughout our paper. In Section 3 we will find the renormalization order of the sequence (\ref{snx}), while Section 4 contains the result on the convergence of the ``bracket'' (\ref{brac}). In Section 5 we will prove the limit theorem in distribution for $S_{n}(0)$  and in Section 6 we will discuss the stable convergence of this sequence.

\section{Preliminaries}

Here we describe the elements from stochastic analysis that we will need in the paper. Consider ${\mathcal{H}}$ a real separable Hilbert space and $(B (\varphi), \varphi\in{\mathcal{H}})$ an isonormal Gaussian process on a probability space $(\Omega, {\cal{A}}, P)$, that is a centered Gaussian family of random variables such that $\mathbf{E}\left( B(\varphi) B(\psi) \right)  = \langle\varphi, \psi\rangle_{{\mathcal{H}}}$. Denote by $I_{n}$ the multiple stochastic integral with respect to
$B$ (see \cite{N}). This $I_{n}$ is actually an isometry between the Hilbert space ${\mathcal{H}}^{\odot n}$(symmetric tensor product) equipped with the scaled norm $\frac{1}{\sqrt{n!}}\Vert\cdot\Vert_{{\mathcal{H}}^{\otimes n}}$ and the Wiener chaos of order $n$ which is defined as the closed linear span of the random variables $H_{n}(B(\varphi))$ where $\varphi\in{\mathcal{H}}, \Vert\varphi\Vert_{{\mathcal{H}}}=1$ and $H_{n}$ is the Hermite polynomial of degree $n\geq 1$
\begin{equation*}
H_{n}(x)=\frac{(-1)^{n}}{n!} \exp \left( \frac{x^{2}}{2} \right)
\frac{d^{n}}{dx^{n}}\left( \exp \left( -\frac{x^{2}}{2}\right)
\right), \hskip0.5cm x\in \mathbb{R}.
\end{equation*}
The isometry of multiple integrals can be written as: for $m,n$ positive integers,
\begin{eqnarray}
\mathbf{E}\left(I_{n}(f) I_{m}(g) \right) &=& n! \langle f,g\rangle _{{\mathcal{H}}^{\otimes n}}\quad \mbox{if } m=n,\nonumber \\
\mathbf{E}\left(I_{n}(f) I_{m}(g) \right) &= & 0\quad \mbox{if } m\not=n.\label{iso}
\end{eqnarray}
It also holds that
\begin{equation*}
I_{n}(f) = I_{n}\big( \tilde{f}\big)
\end{equation*}
where $\tilde{f} $ denotes the symmetrization of $f$ defined by $\tilde{f}%
(x_{1}, \ldots , x_{x}) =\frac{1}{n!} \sum_{\sigma \in {\cal S}_{n}}
f(x_{\sigma (1) }, \ldots , x_{\sigma (n) } ) $.
\\\\
We recall that any square integrable random variable which is measurable with respect to the $\sigma$-algebra generated by $B$ can be expanded into an orthogonal sum of multiple stochastic integrals
\begin{equation}
\label{sum1} F=\sum_{n\geq0}I_{n}(f_{n})
\end{equation}
where $f_{n}\in{\mathcal{H}}^{\odot n}$ are (uniquely determined)
symmetric functions and $I_{0}(f_{0})=\mathbf{E}\left[  F\right]$.
\\\\
Let $L$ be the Ornstein-Uhlenbeck operator
\begin{equation*}
LF=-\sum_{n\geq 0} nI_{n}(f_{n})
\end{equation*}
if $F$ is given by (\ref{sum1}).
\\\\
For $p>1$ and $\alpha \in \mathbb{R}$ we introduce the Sobolev-Watanabe space $\mathbb{D}^{\alpha ,p }$  as the closure of
the set of polynomial random variables with respect to the norm
\begin{equation*}
\Vert F\Vert _{\alpha , p} =\Vert (I -L) ^{\frac{\alpha }{2}} \Vert_{L^{p} (\Omega )}
\end{equation*}
where $I$ represents the identity. We denote by $D$  the Malliavin  derivative operator that acts on smooth functions of the form $F=g(B(\varphi _{1}), \ldots , B(\varphi_{n}))$ ($g$ is a smooth function with compact support and $\varphi_{i} \in {{\cal{H}}}$)
\begin{equation*}
DF=\sum_{i=1}^{n}\frac{\partial g}{\partial x_{i}}(B(\varphi _{1}), \ldots , B(\varphi_{n}))\varphi_{i}.
\end{equation*}
The operator $D$ is continuous from $\mathbb{D} ^{\alpha , p} $ into $\mathbb{D} ^{\alpha -1, p} \left( {\cal{H}}\right).$ The adjoint of $D$ is denoted by $\delta $ and is called the divergence (or Skorohod) integral. It is a continuous operator from $\mathbb{D}^{\alpha, p } \left( {\cal{H}}\right)$ into $\mathbb{D} ^{\alpha -1,p}$.  We have the following duality relationship between $D$ and $\delta$
\begin{equation}
\label{dua}
\mathbf{E} (F\delta (u))= \mathbf{E}\langle DF, u\rangle _{{\cal{H}}} \mbox{ for every $F$ smooth.}
\end{equation}
For adapted integrands, the divergence integral coincides with the classical It\^o integral. We will use the notation
\begin{equation*}
\delta (u) =\int_{0}^{T} u_{s} dB_{s}.
\end{equation*}
Let  $u$ be  a stochastic process having the chaotic decomposition $u_{s}=\sum _{n\geq 0} I_{n}(f_{n}(\cdot ,s))$
where $f_{n}(\cdot, s)\in {\cal{H}}^{\otimes n}$ for every $s$. One can prove that $u \in {\rm Dom} \ \delta$ if and only if $\tilde f_n \in {\cal{H}}^{\otimes (n+1)}$ for every $n \geq 0$, and $\sum_{n=0}^{\infty}I_{n+1}(\tilde f_n)$ converges in $L^2(\Omega)$.
In this case,
$$\delta(u)=\sum_{n=0}^{\infty}I_{n+1}(\tilde f_n) \quad \mbox{and}
\quad \mathbf{E}|\delta(u)|^2=\sum_{n=0}^{\infty}(n+1)! \ \|\tilde
f_n\|_{{\cal{H}}^{\otimes (n+1)}}^{2}.$$
In our  work we will mainly consider divergence integrals with respect to a fractional Brownian motion. The fractional Brownian motion $(B^{H}_{t}) _{t\in [0,T]}$ with  Hurst parameter $H\in (0, 1)$ is a centered Gaussian process starting from zero with covariance function $$R^{H}(t,s):= \frac{1}{2} \left( t^{2H}+s^{2H} -\vert t-s\vert ^{2H}\right), \hskip0.5cm s,t \in [0,T].$$ In this case the space ${\cal{H}_{H}}$ is the canonical Hilbert space of the fractional Brownian motion which is defined as the closure of the linear space generated by the indicator functions $\{ 1_{[0,t] }, t\in [0,T]\}$ with respect to the scalar product
\begin{equation*}
\langle 1_{[0,t]}, 1_{[0,s]} \rangle _{{\cal{H}}_{H}} = R^{H}(t,s), \hskip0.5cm s,t\in [0,T].
\end{equation*}

\section{Renormalization of the sequence  $S_{n}$}
As we mentioned in the introduction, we will assume throughout the paper that $x=0$ in (\ref{snx}), then
\begin{eqnarray}\label{sn}
S_{n}:=S_{n}(0) = \sum_{i=0}^{n-1}K(n^{\alpha}B_{i}^{H_{1}}) (B_{i+1}^{H_{2}} - B_{i}^{H_{2}}).
\end{eqnarray}
We compute in this part the $L^{2}$-norm of $S_{n}$ in order to renormalize it. We have
\begin{eqnarray}
\mathbf{E}\left(S_{n}^{2}\right) & = & \mathbf{E}\left(\sum_{i,j = 0}^{n-1} K(n^{\alpha}B_{i}^{H_{1}}) K(n^{\alpha}B_{j}^{H_{1}}) (B_{i+1}^{H_{2}} - B_{i}^{H_{2}})(B_{j+1}^{H_{2}} - B_{j}^{H_{2}})\right) \nonumber
\\
& = & \mathbf{E}\left(\sum_{i = 0}^{n-1} K^{2}(n^{\alpha}B_{i}^{H_{1}}) (B_{i+1}^{H_{2}} - B_{i}^{H_{2}})^{2}\right) \nonumber
\\
& & + \mathbf{E}\left(\sum_{i \neq j}^{n-1} K(n^{\alpha}B_{i}^{H_{1}}) K(n^{\alpha}B_{j}^{H_{1}}) (B_{i+1}^{H_{2}} - B_{i}^{H_{2}}) (B_{j+1}^{H_{2}} - B_{j}^{H_{2}})\right) \nonumber
\\
& = & T' + T''. \nonumber
\end{eqnarray}
The summand $T'$ will be called the ``diagonal'' term while the summand $T''$ will be called ``the non-diagonal'' term. We will analyze each of them separately. Concerning $T'$ we have

\begin{lemma}\label{l1}
As $n \rightarrow + \infty$,
\begin{eqnarray}
n^{\alpha + H_{1} - 1}T' \underset{n \rightarrow +\infty}{\longrightarrow} C_{1} = \frac{1}{2\pi \sqrt{2}(1-H_{1})}.
\end{eqnarray}
\end{lemma}
\begin{dem}
Through the independence of $\left(B_{t}^{H_{1}}\right)_{t \geq 0}$ and $\left(B_{t}^{H_{2}}\right)_{t \geq 0}$,
\begin{eqnarray}
T' = \sum_{i = 0}^{n-1} \mathbf{E}\left(K^{2}(n^{\alpha}B_{i}^{H_{1}})\right) \mathbf{E}\left((B_{i+1}^{H_{2}} - B_{i}^{H_{2}})^{2}\right). \nonumber
\end{eqnarray}
Since $\mathbf{E}\left((B_{i+1}^{H_{2}} - B_{i}^{H_{2}})^{2}\right) = 1$,
\begin{eqnarray}
T' = \sum_{i = 0}^{n-1} \mathbf{E}\left(K^{2}(n^{\alpha}B_{i}^{H_{1}})\right) = \sum_{i = 0}^{n-1} \mathbf{E}\left(\frac{1}{2\pi}e^{-n^{2\alpha}i^{2H_{1}}Z^{2}}\right) \nonumber
\end{eqnarray}
where $Z$ is a standard normal random variable. Recall that, if $Z$ is a standard normal random variable, and if $1+2c>0$
\begin{equation}
\label{exp}
\mathbf{E}\left( e^{-cZ^{2}} \right) =\frac{1}{\sqrt{1+2c}}
\end{equation}
consequently,
\begin{eqnarray}
T' = \sum_{i = 0}^{n-1} \frac{1}{2\pi \sqrt{1 + 2 n^{2\alpha}i^{2H_{1}}}}. \nonumber
\end{eqnarray}
As $n \to +\infty$, $T'$ behaves as such
\begin{eqnarray}
\sum_{i = 0}^{n-1} \frac{1}{2\pi \sqrt{1 + 2 n^{2\alpha}i^{2H_{1}}}} & \sim & \frac{n^{-\alpha}}{2\pi \sqrt{2}}\sum_{i = 0}^{n-1}i^{-H_{1}} \sim \frac{n^{-\alpha - H_{1} +1}}{2\pi \sqrt{2}}\frac{1}{n}\sum_{i = 0}^{n-1}\left(\frac{i}{n}\right)^{-H_{1}}\nonumber
\\
& \sim & \frac{n^{-\alpha - H_{1} +1}}{2\pi \sqrt{2}} \int_{0}^{1}x^{- H_{1}}dx = \frac{n^{-\alpha - H_{1} +1}}{2\pi \sqrt{2}(1-H_{1})}. \nonumber
\end{eqnarray}
The sign ``$\sim$'' means that the left-hand side and the right-hand side have the same limit as $n \to +\infty$. We will use this notation throughout the paper.
\qed
\end{dem}

\vskip0.3cm

\noindent We will now compute the term $T''$. To do so, we will need the following Lemma (lemma 3.1 p. 122 in \cite{YanLiuYang}).
\begin{lemma}
\label{majoreDet}
For every $s,r \in \left[0,T\right]$, $s \geq r$ and $0<H<1$ we have
\begin{eqnarray}
s^{2H}r^{2H} - \mu^{2} \geq \tau(s-r)^{2H}r^{2H}
\end{eqnarray}
where $\mu = \mathbf{E}(B_{s}^{H}B_{r}^{H})$ and $\tau > 0$ is a constant.
\end{lemma}
Concerning the non-diagonal term of $\mathbf{E}\left(S_{n}^{2}\right)$ the following holds
\begin{lemma}\label{l2}
Suppose that
\begin{equation}
\label{cond}
\alpha -4H_{2}+H_{1}+2>0.
\end{equation}
Then, as $n \rightarrow + \infty$,
\begin{eqnarray}
n^{\alpha +H_{1}-1}T'' \underset{n \rightarrow +\infty}{\longrightarrow}0.
\end{eqnarray}
\end{lemma}
\begin{dem}
Using again the independence of  $\left(B_{t}^{H_{1}}\right)_{t \geq 0}$ and $\left(B_{t}^{H_{2}}\right)_{t \geq 0}$
\begin{eqnarray}
T'' & = & \sum_{i \neq j}^{n-1} \mathbf{E}\left(K(n^{\alpha}B_{i}^{H_{1}}) K(n^{\alpha}B_{j}^{H_{1}})\right)\mathbf{E}\left((B_{i+1}^{H_{2}} - B_{i}^{H_{2}}) (B_{j+1}^{H_{2}} - B_{j}^{H_{2}})\right) \nonumber
\\
& = & \frac{1}{2}\sum_{i \neq j}^{n-1} \mathbf{E}\left(K(n^{\alpha}B_{i}^{H_{1}}) K(n^{\alpha}B_{j}^{H_{1}})\right)f_{H_{2}}(i,j) \nonumber
\end{eqnarray}
where
\begin{equation}
\label{notationCovariancefBm}
f_{H_{2}}(i,j)=\frac{1}{2}\left[\left|i - j + 1\right|^{2H_{2}} + \left|i - j - 1\right|^{2H_{2}} -2\left|i - j\right|^{2H_{2}}\right].
\end{equation}
We need to evaluate the expectation  $\mathbf{E}\left(K(n^{\alpha}B_{i}^{H_{1}}) K(n^{\alpha}B_{j}^{H_{1}})\right)$. Let $\Gamma=\begin{pmatrix}
i^{2H_{1}} & R(i,j) \\
R(i,j) & j^{2H_{1}}
\end{pmatrix}$ be the covariance matrix of  $\left(B_{i}^{H_{1}},B_{j}^{H_{1}}\right)$. We have $\left|\Gamma\right| = (ij)^{2H_{1}} - R^{2}(i,j)$ and $\Gamma^{-1}=\frac{1}{\left|\Gamma\right|}\begin{pmatrix}
j^{2H_{1}} & -R(i,j) \\
-R(i,j) & i^{2H_{1}}
\end{pmatrix}$. The density of $\left(B_{i}^{H_{1}},B_{j}^{H_{1}}\right)$ is then
\begin{eqnarray}
f(x,y) = \frac{1}{2\pi \sqrt{\left|\Gamma\right|}}e^{-\frac{1}{2\left|\Gamma\right|}(j^{2H_{1}}x^{2} - 2R(i,j)xy + i^{2H_{1}}y^{2})}.
\end{eqnarray}
We obtain
\begin{eqnarray}
&& \mathbf{E}\left(K(n^{\alpha}B_{i}^{H_{1}}) K(n^{\alpha}B_{j}^{H_{1}})\right)=  \frac{1}{(2\pi)^{2} \sqrt{\left|\Gamma\right|}}\int_{\mathbb{R}^{2}}e^{-\frac{n^{2\alpha}x^{2}}{2}}e^{-\frac{n^{2\alpha}y^{2}}{2}}e^{-\frac{1}{2\left|\Gamma\right|}(j^{2H_{1}}x^{2} - 2R(i,j)xy + i^{2H_{1}}y^{2})}dxdy \nonumber
\\
&=& \frac{1}{(2\pi)^{2} \sqrt{\left|\Gamma\right|}} \int_{\mathbb{R}}e^{-\frac{n^{2\alpha}y^{2}}{2}}e^{-\frac{i^{2H_{1}}y^{2}}{2\left|\Gamma\right|}} \int_{\mathbb{R}}e^{-\frac{n^{2\alpha}x^{2}}{2}}e^{-\frac{1}{2\left|\Gamma\right|}(j^{2H_{1}}x^{2} - 2R(i,j)xy )}dxdy \nonumber
\\
&=& \frac{1}{(2\pi)^{2} \sqrt{\left|\Gamma\right|}} \int_{\mathbb{R}}e^{-\frac{y^{2}}{2}\left[n^{2\alpha} + \frac{i^{2H_{1}}}{\left|\Gamma\right|} \right]} \int_{\mathbb{R}}e^{-\frac{1}{2}\left[x^{2}\left(n^{2\alpha} + \frac{j^{2H_{1}}}{\left|\Gamma\right|}\right) -\frac{2R(i,j)}{\left|\Gamma\right|}xy\right]}dxdy \nonumber
\\
&=& \frac{1}{(2\pi)^{2} \sqrt{\left|\Gamma\right|}} \int_{\mathbb{R}}e^{-\frac{y^{2}}{2}\left[n^{2\alpha} + \frac{i^{2H_{1}}}{\left|\Gamma\right|} \right]} \int_{\mathbb{R}}e^{-\frac{\left(n^{2\alpha} + \frac{j^{2H_{1}}}{\left|\Gamma\right|}\right)}{2}\left[x^{2} - \frac{2R(i,j)}{n^{2\alpha}\left|\Gamma\right| + j^{2H_{1}}}xy\right]}dxdy \nonumber
\\
&=& \frac{1}{(2\pi)^{2} \sqrt{\left|\Gamma\right|}} \int_{\mathbb{R}}e^{-\frac{y^{2}}{2}\left[n^{2\alpha} + \frac{i^{2H_{1}}}{\left|\Gamma\right|} \right]} \int_{\mathbb{R}}e^{-\frac{\left(n^{2\alpha} + \frac{j^{2H_{1}}}{\left|\Gamma\right|}\right)}{2}\left[\left(x - \frac{R(i,j)}{n^{2\alpha}\left|\Gamma\right| + j^{2H_{1}}}y\right)^{2} - \frac{R^{2}(i,j)}{(n^{2\alpha}\left|\Gamma\right| + j^{2H_{1}})^{2}}y^{2}\right]}dxdy \nonumber
\\
&=& \frac{1}{(2\pi)^{2} \sqrt{\left|\Gamma\right|}} \int_{\mathbb{R}}e^{-\frac{y^{2}}{2}\left[n^{2\alpha} + \frac{i^{2H_{1}}}{\left|\Gamma\right|} \right]}e^{-\frac{\left(n^{2\alpha} + \frac{j^{2H_{1}}}{\left|\Gamma\right|}\right)}{2}\frac{R^{2}(i,j)y^{2}}{\left(n^{2\alpha} + \frac{j^{2H_{1}}}{\left|\Gamma\right|}\right)^{2}\left|\Gamma\right|^{2}}} \int_{\mathbb{R}}e^{-\frac{\left(n^{2\alpha} + \frac{j^{2H_{1}}}{\left|\Gamma\right|}\right)}{2}\left(x - \frac{R(i,j)}{n^{2\alpha}\left|\Gamma\right| + j^{2H_{1}}}y\right)^{2}}dxdy \nonumber
\\
&=& \frac{1}{(2\pi)^{\frac{3}{2}} \sqrt{\left|\Gamma\right|}}\frac{\sqrt{\left|\Gamma\right|}}{\sqrt{n^{2\alpha}\left|\Gamma\right| + j^{2H_{1}}}} \int_{\mathbb{R}}e^{-\frac{1}{2}y^{2}\left[\frac{(n^{2\alpha}\left|\Gamma\right| + i^{2H_{1}})(n^{2\alpha}\left|\Gamma\right| + j^{2H_{1}}) - R^{2}(i,j)}{\left|\Gamma\right|(n^{2\alpha}\left|\Gamma\right| + j^{2H_{1}})}\right]}dy.   \nonumber
\end{eqnarray}
Thus
\begin{eqnarray}
\mathbf{E}\left(K(n^{\alpha}B_{i}^{H_{1}}) K(n^{\alpha}B_{j}^{H_{1}})\right) &=& \frac{1}{2\pi\sqrt{n^{2\alpha}\left|\Gamma\right| + j^{2H_{1}}}}\frac{\sqrt{\left|\Gamma\right|}\sqrt{(n^{2\alpha}\left|\Gamma\right| + j^{2H_{1}})}}{\sqrt{(n^{2\alpha}\left|\Gamma\right| + i^{2H_{1}})(n^{2\alpha}\left|\Gamma\right| + j^{2H_{1}}) - R^{2}(i,j)}} \nonumber
\\
& = & \frac{\sqrt{\left|\Gamma\right|}}{2\pi \sqrt{(n^{2\alpha}\left|\Gamma\right| + i^{2H_{1}})(n^{2\alpha}\left|\Gamma\right| + j^{2H_{1}}) - R^{2}(i,j)}} \nonumber
\\
& = & \frac{1}{2\pi \sqrt{n^{4\alpha}\left|\Gamma\right| + n^{2\alpha}j^{2H_{1}} + n^{2\alpha}i^{2H_{1}} + 1}}. \nonumber
\end{eqnarray}
Suppose that $i>j$. We use Lemma \ref{majoreDet} to bound $\left|\Gamma\right| = i^{2H_{1}}j^{2H_{1}} - R^{2}(i,j)$ from below. Therefore
\begin{eqnarray}
\mathbf{E}\left(K(n^{\alpha}B_{i}^{H_{1}}) K(n^{\alpha}B_{j}^{H_{1}})\right) & \leq & \frac{1}{2\pi \sqrt{n^{4\alpha}\tau (i-j)^{2H_{1}}j^{2H_{1}} + n^{2\alpha}(i^{2H_{1}} + j^{2H_{1}})}}. \nonumber
\end{eqnarray}
 Since $a^{2}+b^{2}\geq 2ab$ with $a^{2} = n^{4\alpha}\tau (i-j)^{2H_{1}}j^{2H_{1}}$ and $b^{2} = n^{2\alpha}(i^{2H_{1}} + j^{2H_{1}})$
\begin{eqnarray}
\mathbf{E}\left(K(n^{\alpha}B_{i}^{H_{1}}) K(n^{\alpha}B_{j}^{H_{1}})\right) & \leq & \frac{1}{2\pi \sqrt{2\sqrt{\tau}n^{2\alpha} (i-j)^{H_{1}}j^{H_{1}}\sqrt{n^{2\alpha}(i^{2H_{1}} + j^{2H_{1}})}}} \nonumber
\end{eqnarray}
and using the same inequality as above for $a^{2} = i^{2H_{1}}$ and $b^{2} = j^{2H_{1}}$
\begin{eqnarray}
\label{eqMajT2}
\mathbf{E}\left(K(n^{\alpha}B_{i}^{H_{1}}) K(n^{\alpha}B_{j}^{H_{1}})\right) & \leq & \frac{n^{-\frac{3\alpha}{2}}}{2\pi\sqrt{2}\tau^{\frac{1}{4}}(i-j)^{\frac{H_{1}}{2}}j^{\frac{3H_{1}}{4}}i^{\frac{H_{1}}{4}}}.
\end{eqnarray}
Since $f_{H_{2}}(i,j)$ behaves as $H_{2}(2H_{2}-1)\vert i-j\vert ^{2H_{2}-2}$ when $i-j\to \infty$, we can assert that
\begin{eqnarray}
T'' & \sim & \frac{H_{2}(2H_{2} -1)}{2}\sum_{i \neq j}^{n-1} \mathbf{E}\left(K(n^{\alpha}B_{i}^{H_{1}}) K(n^{\alpha}B_{j}^{H_{1}})\right)\left|i - j\right|^{2H_{2} -2}. \nonumber
\end{eqnarray}
Using (\ref{eqMajT2}), we can write
\begin{eqnarray}
\sum_{i \neq j}^{n-1} \mathbf{E}\left(K(n^{\alpha}B_{i}^{H_{1}}) K(n^{\alpha}B_{j}^{H_{1}})\right)\left|i - j\right|^{2H_{2} -2} & \lesssim & \sum_{i > j}^{n-1} \frac{n^{-\frac{3\alpha}{2}}}{2\pi\sqrt{2}\tau^{\frac{1}{4}}(i-j)^{\frac{H_{1}}{2}}j^{\frac{3H_{1}}{4}}i^{\frac{H_{1}}{4}}}\left|i - j\right|^{2H_{2} -2}\nonumber
\end{eqnarray}
and consequently
\begin{eqnarray}
T'' & \lesssim & \frac{H_{2}(2H_{2} -1)}{4\pi \sqrt{2} \tau^{\frac{1}{4}}}n^{-\frac{3\alpha}{2}}n^{2H_{2}-\frac{H_{1}}{2} - 2}n^{-\frac{3H_{1}}{4}}n^{-\frac{H_{1}}{4}}n^{2}\underbrace{\frac{1}{n^{2}}\sum_{i > j}^{n-1}\frac{\left(\frac{i-j}{n}\right)^{2H_{2}-\frac{H_{1}}{2} - 2}}{\left(\frac{j}{n}\right)^{\frac{3H_{1}}{4}}}\left(\frac{j}{n}\right)^{\frac{H_{1}}{4}}}_{\underset{n \rightarrow +\infty}{\longrightarrow} C(H_{1},H_{2}) > 0} \nonumber
\\
& \lesssim & \frac{H_{2}(2H_{2} -1)C(H_{1},H_{2})}{4\pi \sqrt{2} \tau^{\frac{1}{4}}}n^{-\frac{3\alpha}{2} +2H_{2}-\frac{3H_{1}}{2}}.
\end{eqnarray}
It follows that under condition (\ref{cond}) $ n^{\alpha+H_{1}-1}T''$
converges to zero as $n\to \infty$. \qed
\end{dem}

\vskip0.3cm

\noindent As a consequence of Lemmas \ref{l1} and \ref{l2} we obtain the following $L^{2}$- norm estimate for $S_{n}$.
\begin{prop}\label{prop1}
Suppose that condition (\ref{cond}) holds. Then, as $n\to \infty$
\begin{equation*}
n^{\alpha+H_{1}-1}\mathbf{E}\left(S_{n}^{2}\right) \to C_{1} = \frac{1}{2\pi \sqrt{2}(1-H_{1})}.
\end{equation*}
\end{prop}
The condition (\ref{cond}) will be discussed more thoroughly later (Remark 1, Section 5).

\section{The limit in distribution of $\langle S\rangle _{n}$}

Proposition \ref{prop1} implies that the diagonal part of $S_{n}^{2}$ is dominant in relation to the non-diagonal part, in the sense that this diagonal part is responsable for the renormalization order of $S_{n}^{2}$ which is $n^{\alpha +H_{1}-1}$. As a consequence we need to study the limit distribution of $n^{\alpha +H_{1}-1}\langle S\rangle _{n}= n^{\alpha +H_{1}-1}\sum _{i=0}^{n-1}K^{2}(n^{\alpha } B^{H_{1}}_{i})$. Using the self-similarity property  of the fractional Brownian motion we have
\begin{equation*}
n^{\alpha +H_{1}-1}\sum _{i=0}^{n-1}K^{2}(n^{\alpha } B^{H_{1}}_{i})=n^{\alpha +H_{1}-1}\sum_{i=0}^{n-1} K^{2} (n^{\alpha +
H_{1}}B_{\frac{i}{n}}^{H_{1}}).
\end{equation*}
The limit of the above sequence is linked to the local time of the fractional Brownian motion $B^{H_{1}}$.  For any $t\geq 0$ and $x\in \mathbb{R}$  we define $L^{H_{1}}(t,x)$ as the density of the occupation measure (see \cite{Be}, \cite{GH})
\begin{equation*}
\mu _{t}(A)= \int_{0}^{t} 1_{A}(B^{H_{1}}_{s})ds, \hskip0.5cm A\in {\cal{B}}(\mathbb{R}).
\end{equation*}
The local time $L^{H_{1}}(t,x)$ satisfies the occupation time formula
\begin{equation}
\label{oc}
\int_{0}^{t} f(B^{H_{1}}_{s})ds = \int_{\mathbb{R}}L^{H_{1}}(t,x) f(x)dx
\end{equation}
for any measurable function $f$. The local time is H\"older continuous with respect to $t$ and with respect to $x$ (for the sake of completeness  $L^{H_{1}}(t,x)$ has H\"older continuous paths of order $\delta <1-H$ in time and of order $\gamma <\frac{1-H}{2H} $  in the space variable  (see Table 2 in \cite{GH})).  Moreover, it admits a bicontinuous version with respect to $(t,x)$.
\\\\
\noindent Below, we give an important convergence result that will be necessary in proving the main result of this section.
\begin{prop}\label{prop2}
The following convergence in distribution result holds
\begin{eqnarray}
\label{riemannSumConv}
n^{\alpha + H_{1}}\left( \frac{1}{n}\sum_{i=0}^{n-1}K^{2}(n^{\alpha +
H_{1}}B_{\frac{i}{n}}^{H_{1}}) - \int_{0}^{1}K^{2}(n^{\alpha +
H_{1}}B_{s}^{H_{1}})ds \right) \underset{n \rightarrow
+\infty}{\longrightarrow} 0.
\end{eqnarray}
\end{prop}
\begin{dem}
\noindent Fix $\varepsilon>0$. Let $p_{\varepsilon}(x)$ be the
Gaussian kernel with variance $\varepsilon >0$ defined by $p_{\varepsilon} (x) = \frac{1}{\sqrt{2\pi\varepsilon}}e^{-\frac{x^{2}}{2\eps}}$.  Note that for every $s\geq 0$
\begin{eqnarray}
\label{kern}
\sqrt{\pi}n^{\alpha +H_{1}} K^{2}(n^{\alpha + H_{1}}B_{s}^{H_{1}})=\frac{1}{2} p_{\frac{1}{2n^{2(\alpha + H_{1})}}}(B^{H_{1} }_{s}).
\end{eqnarray}
Using (\ref{kern}), we can write the left-hand side of (\ref{riemannSumConv}) as
\begin{eqnarray*}
&&\sqrt{\pi} n^{\alpha + H_{1}}\left( \int_{0}^{1}K^{2}(n^{\alpha + H_{1}}B_{s}^{H_{1}})ds -\frac{1}{n}\sum_{i=0}^{n-1}K^{2}(n^{\alpha + H_{1}}B_{\frac{i}{n}}^{H_{1}})\right) \\
&=&\frac{1}{2}\sum_{i=0}^{n-1}\int_{\frac{i}{n}}^{\frac{i+1}{n}} \left( p_{\frac{1}{2} n^{-2(\alpha +H_{1})}}(B^{H_{1}}_{s}) -p_{\frac{1}{2} n^{-2(\alpha +H_{1})}}(B^{H_{1}}_{\frac{i}{n}})\right)ds\\
&=&\frac{1}{2}\sum_{i=0}^{n-1}\int_{\frac{i}{n}}^{\frac{i+1}{n}} \left( p_{\frac{1}{2} n^{-2(\alpha +H_{1})}}(B^{H_{1}}_{s})-p_{\varepsilon} (B^{H_{1}}_{s}) \right)ds \\
&&+\frac{1}{2}\sum_{i=0}^{n-1}\int_{\frac{i}{n}}^{\frac{i+1}{n}}\left( p_{\varepsilon} (B^{H_{1}}_{s}) -p_{\varepsilon}(B^{H_{1}}_{\frac{i}{n}})\right) ds \\
&&+\frac{1}{2}\sum_{i=0}^{n-1}\int_{\frac{i}{n}}^{\frac{i+1}{n}}\left( p_{\varepsilon}(B^{H_{1}}_{\frac{i}{n}})-p_{\frac{1}{2}n^{-2(\alpha +H_{1})}}(B^{H_{1}}_{\frac{i}{n}})\right)ds: = \frac{1}{2}(a_{n}^{(1)}+ a_{n}^{(2)} + a_{n}^{(3)}).
\end{eqnarray*}
We will now estimate the three terms above and we will show that
each of them converges to zero (in some sense). Let us first handle
the term $a_{n}^{(1)}$. We have
$$a_{n}^{(1)}= \int_{0}^{1} p_{\frac{1}{2}n^{-2(\alpha +H_{1})}} (B^{H_{1}}_{s}) ds -\int_{0}^{1} p_{\varepsilon} (B^{H_{1}}_{s})ds.$$
It follows from \cite{NV} or \cite{Edd} that
\begin{eqnarray}
\int_{0}^{1} p_{\varepsilon} (B^{H_{1}}_{s})ds\to _{\varepsilon \to
0} \int_{0}^{1} \delta _{0} (B^{H_{1}}_{s}) ds =L^{H_{1}}(1,0)
\end{eqnarray}
in $L^{2}(\Omega)$ and almost surely, where $L^{H_{1}}(1,0)$ is the local time of the
fractional Brownian motion. Therefore $a_{n}^{(1)}$ clearly converges to zero as $\varepsilon \to 0$  and $n\to \infty$. The
term $a_{n}^{(2)}$ can be expressed as
\begin{eqnarray}
a_{n}^{(2)}= -\left(\frac{1}{n}\sum_{i=0}^{n-1} p_{\varepsilon}
(B^{H_{1}}_{\frac{i}{n}})-\int_{0}^{1}
p_{\varepsilon}(B^{H_{1}}_{s}) ds\right)
\end{eqnarray}
and for every $\varepsilon >0$ it converges almost surely to zero as
$n\to \infty$ using the Riemann sum convergence. Let us now handle
the term $a_{n}^{(3)}$ given by
\begin{equation}\label{a3}a_{n}^{(3)}= \frac{1}{n} \sum_{i=0}^{n-1} \left( p_{\varepsilon} (B^{H_{1}}_{\frac{i}{n}})- p_{\frac{1}{2} n^{-2(\alpha +H_{1})}}(B^{H_{1}}_{\frac{i}{n}})\right).
  \end{equation}
  We will treat this term by using the chaos decomposition of the Gaussian kernel applied to random variables in the first Wiener chaos. Recall that (see \cite{CNT}, \cite{HO}, \cite{Imk1}, \cite{NV2}) for every $\varphi \in {\cal{H}}_{H_{1}}$ (${\cal{H}}_{H_{1}}$ is the canonical Hilbert space associated with the Gaussian process $B^{H_{1}}$),
\begin{eqnarray}
p_{\varepsilon}(B^{H_{1}}(\varphi))=\sum_{m\geq 0}C_{m}I_{2m} \left(
\varphi ^{\otimes 2m} \right) \frac{1}{\left(\Vert \varphi
\Vert^{2}_{{\cal{H}}_{1}}+ \varepsilon\right) ^{m+\frac{1}{2}} }
\end{eqnarray}
where $C_{m} = \frac{(-1)^{m}}{\sqrt{2\pi}2^{m}m!}$.
\\\\
\noindent Using this chaos decomposition, we can write $p_{\varepsilon}
(B^{H_{1}}_{\frac{i}{n}})- p_{\frac{1}{2} n^{-2(\alpha
+H_{1})}}(B^{H_{1}}_{\frac{i}{n}})$ as
\begin{eqnarray*}
p_{\varepsilon} (B^{H_{1}}_{\frac{i}{n}})- p_{\frac{1}{2}
n^{-2(\alpha +H_{1})}}(B^{H_{1}}_{\frac{i}{n}})&=&\sum_{m\geq
0}C_{m}I_{2m} \left( 1_{[0, \frac{i}{n}]}^{\otimes 2m}\right) \left(
\frac{1}{ \left( \left( \frac{i}{n}\right) ^{2H_{1}} +
\varepsilon\right) ^{m+\frac{1}{2}}}-\frac{1}{ \left( \left(
\frac{i}{n}\right) ^{2H_{1}} + \frac{1}{2}n^{-2(\alpha +
H_{1})}\right) ^{m+\frac{1}{2}}}\right)
\\
&=& \sum_{m\geq 0}C_{m}I_{2m} \left( 1_{[0, \frac{i}{n}]}^{\otimes
2m}\right)\left(\frac{i}{n}\right)^{-2H_{1}\left(m+\frac{1}{2}\right)}
d_{i,\varepsilon,n,m}
\end{eqnarray*}
where
\begin{eqnarray*}
d_{i,\varepsilon,n,m} = \left(\left(\frac{\left( \frac{i}{n}\right)
^{2H_{1}}}{ \left( \left( \frac{i}{n}\right) ^{2H_{1}} +
\varepsilon\right)}\right)^{m+\frac{1}{2}} - \left(\frac{\left(
\frac{i}{n}\right) ^{2H_{1}}}{ \left( \left( \frac{i}{n}\right)
^{2H_{1}} + \frac{1}{2}n^{-2(\alpha +
H_{1})}\right)}\right)^{m+\frac{1}{2}}\right).
\end{eqnarray*}
We will show that $a_{n}^{(3)}$ converges to zero in $L^{2}(\Omega)$
as $n\to \infty$ and $\varepsilon \to 0$. From (\ref{a3}) one can easily see that the diagonal part of $a_{n}^{(3)}$ converges to zero. We can also see, from the expression of $a_{n}^{(3)}$, that the summands with $j=0$ vanish. Then, by using the
orthogonality of multiple stochastic integrals(\cite{N}), we obtain
\begin{eqnarray*}
\mathbf{E}(a_{n}^{(3)})^{2} \sim \frac{1}{n^{2}} \sum_{m\geq 0}
C_{m}^{2}(2m)!\sum_{i,j\geq 1; i\not=j}^{n-1} \langle 1_{[0, \frac{i}{n} ]} ,
1_{[0, \frac{j}{n} ]}\rangle _{{\cal{H}}_{1}}^{2m}
\left(\frac{i}{n}\right) ^{-2H_{1} (m+\frac{1}{2})
}\left(\frac{j}{n}\right) ^{-2H_{1} (m+\frac{1}{2})
}d_{i,\varepsilon, n,m } d_{j,\varepsilon, n, m}.
\end{eqnarray*}
We can also write
\begin{eqnarray*}
\mathbf{E}(a_{n}^{(3)})^{2} &\sim& \frac{1}{n^{2}} \sum_{m\geq 0}
C_{m}^{2}(2m)!\sum_{i,j\geq 1, i\not= j}^{n-1}
R_{H_{1}}\left(\frac{i}{n},\frac{j}{n}\right)^{2m}
\left(\frac{i}{n}\right) ^{-2H_{1} (m+\frac{1}{2})
}\left(\frac{j}{n}\right) ^{-2H_{1} (m+\frac{1}{2})
}d_{i,\varepsilon, n,m } d_{j,\varepsilon, n, m}
\\
&:=& \sum_{m\geq 0} C_{m}^{2}(2m)!A_{m}(\varepsilon, n).
\end{eqnarray*}
where
$$A_{m}(\varepsilon, n) =\frac{1}{n^{2}}\sum_{i,j\geq 1; i\not= j}^{n-1} R_{H_{1}}\left( \frac{i}{n}, \frac{j}{n} \right) ^{2m} \left( \frac{i}{n}\right) ^{-2H_{1} (m+\frac{1}{2}) }\left( \frac{j}{n}\right) ^{-2H_{1} (m+\frac{1}{2}) }d_{i,\varepsilon, n,m}d_{j,\varepsilon, n, m}.$$
We can now claim that, for every fixed $m\geq 1$
\begin{equation}
\label{ame} \lim_{\varepsilon \to 0} \lim _{n\to \infty}
A_{m}(\varepsilon,n )=0.
\end{equation}
Indeed, for every $m\geq 0$, we get
\begin{eqnarray*}
|d_{i,\varepsilon, n,m}|&=& \left| \left(\left(\frac{\left( \frac{i}{n}\right)
^{2H_{1}}}{ \left( \left( \frac{i}{n}\right) ^{2H_{1}} +
\varepsilon\right)}\right)^{m+\frac{1}{2}} -1 +1 -\left(\frac{\left(
\frac{i}{n}\right) ^{2H_{1}}}{ \left( \left( \frac{i}{n}\right)
^{2H_{1}} + \frac{1}{2}n^{-2(\alpha +
H_{1})}\right)}\right)^{m+\frac{1}{2}}\right) \right| \\
&\leq & \left| 1-\left(\frac{\left( \frac{i}{n}\right)
^{2H_{1}}}{ \left( \left( \frac{i}{n}\right) ^{2H_{1}} +
\varepsilon\right)}\right)^{m+\frac{1}{2}} \right|+ \left| 1 -\left(\frac{\left(
\frac{i}{n}\right) ^{2H_{1}}}{ \left( \left( \frac{i}{n}\right)
^{2H_{1}} + \frac{1}{2}n^{-2(\alpha +
H_{1})}\right)}\right)^{m+\frac{1}{2}}\right|\\
&\leq & \left| 1-\left(\frac{\left( \frac{i}{n}\right)
^{2H_{1}}}{ \left( \left( \frac{i}{n}\right) ^{2H_{1}} +
\varepsilon\right)}\right)^{m+1} \right|+ \left| 1 -\left(\frac{\left(
\frac{i}{n}\right) ^{2H_{1}}}{ \left( \left( \frac{i}{n}\right)
^{2H_{1}} + \frac{1}{2}n^{-2(\alpha +
H_{1})}\right)}\right)^{m+1}\right|\\
&=&c_{m}\left( \left| \left(\frac{\varepsilon}{ \left( \left( \frac{i}{n}\right) ^{2H_{1}} +
\varepsilon\right)}\right) \right|+
 \left| \left(\frac{n^{-2(\alpha +H_{1})}}{ \left( \left( \frac{i}{n}\right)
^{2H_{1}} + \frac{1}{2}n^{-2(\alpha +
H_{1})}\right)}\right)\right|\right).
 \end{eqnarray*}
 Now, for every $i,n,m$, we have $\underset{\eps \to 0 }{\mbox{lim}}\left| \left(\frac{\varepsilon}{ \left( \left( \frac{i}{n}\right) ^{2H_{1}} +
\varepsilon\right)}\right) \right| =0$ and for every $i\geq 1$,
$$\left| \left(\frac{n^{-2(\alpha +H_{1})}}{ \left( \left( \frac{i}{n}\right)
^{2H_{1}} + \frac{1}{2}n^{-2(\alpha +
H_{1})}\right)}\right)\right| \leq \left| \left(\frac{n^{-2(\alpha +H_{1})}}{ \left( \left( \frac{1}{n}\right)
^{2H_{1}} + \frac{1}{2}n^{-2(\alpha +
H_{1})}\right)}\right)\right|\leq c\frac{n^{2H_{1}}}{n^{2(\alpha +2H_{1})}+n^{2H_{1}}}\underset{n \rightarrow +\infty}{\longrightarrow} 0$$
because $\alpha >0$.
\\\\
\noindent Furthermore, we know that $$\frac{1}{n^{2}}\sum_{i,j=0}^{n-1} R_{H_{1}}\left(
\frac{i}{n}, \frac{j}{n} \right) ^{2m} \left( \frac{i}{n}\right)
^{-2H_{1} (m+\frac{1}{2}) }\left( \frac{j}{n}\right) ^{-2H_{1}
(m+\frac{1}{2}) }$$ converges as $n\to \infty$ to $\int_{0}^{1}\int_{0}^{1}  R(u,v)^{2m}  (uv) ^{-2H_{1}(m+\frac{1}{2})}dudv.$ Since this quantity is finite (\cite{CNT} and \cite{Edd}), it implies  (\ref{ame}).\\\\
\noindent We will now prove that
\begin{equation}
\label{cn} \sum_{m\geq 0} C_{m}^{2}(2m)! \sup_{n,\varepsilon} \left|
A_{m}(\varepsilon ,n) \right| <\infty.
\end{equation}
Relation (\ref{ame}) and (\ref{cn}) will imply the convergence of $a_{n}^{(3)}$ to zero in $L^{2}(\Omega)$.
We need to find an upper bound for the terms $\vert
d_{i,\varepsilon, n,m}\vert $ and $\vert d_{j,\varepsilon,
n,m}\vert$ in order to continue.
\begin{eqnarray*}
d_{i,\varepsilon, n,m} & = & \left(\left(\frac{\left(
\frac{i}{n}\right) ^{2H_{1}}}{ \left( \left( \frac{i}{n}\right)
^{2H_{1}} + \varepsilon\right)}\right)^{m+\frac{1}{2}} -
\left(\frac{\left( \frac{i}{n}\right) ^{2H_{1}}}{ \left( \left(
\frac{i}{n}\right) ^{2H_{1}} + \frac{1}{2}n^{-2(\alpha +
H_{1})}\right)}\right)^{m+\frac{1}{2}}\right)
\\
& = & \left(\left(\frac{1}{ \left( 1 + \varepsilon
n^{2H}i^{-2H}\right)}\right)^{m+\frac{1}{2}} - \left(\frac{1}{
\left( 1 +
\frac{1}{2}n^{-2\alpha}i^{-2H}\right)}\right)^{m+\frac{1}{2}}\right).
\end{eqnarray*}
One can note that
\begin{eqnarray*}
0 \leq \left(\frac{1}{ \left( 1 + \varepsilon
n^{2H}i^{-2H}\right)}\right)^{m+\frac{1}{2}} \leq 1 \mbox{\   \   and\   \   } 0 \leq \left(\frac{1}{ \left( 1 +
\frac{1}{2}n^{-2\alpha}i^{-2H}\right)}\right)^{m+\frac{1}{2}} \leq 1
\end{eqnarray*}
because $\varepsilon n^{2H}i^{-2H} >0$. From the above inequalities, we can deduce that
\begin{eqnarray*}
-1 \leq \left(\frac{1}{ \left( 1 + \varepsilon
n^{2H}i^{-2H}\right)}\right)^{m+\frac{1}{2}} - \left(\frac{1}{
\left( 1 +
\frac{1}{2}n^{-2\alpha}i^{-2H}\right)}\right)^{m+\frac{1}{2}} \leq 1
\end{eqnarray*}
and finally,
\begin{eqnarray*}
\left| d_{i,\varepsilon, n,m}\right|  \leq  1 \mbox{\   \   and\   \   } \left| d_{j,\varepsilon, n,m}\right|  \leq  1.
\end{eqnarray*}
By bounding from above the terms $\vert d_{i,\varepsilon, n,m}\vert $ and
$\vert d_{j,\varepsilon, n,m}\vert$ by 1 in $\sum_{m\geq 0}
C_{m}^{2}(2m)! \sup_{n,\varepsilon} \left| A_{m}(\varepsilon ,n)
\right|$ we obtain that
\begin{eqnarray*}
\sum_{m\geq 0} C_{m}^{2}(2m)! \sup_{n,\varepsilon} \left|
A_{m}(\varepsilon ,n) \right| & \leq &
\sum_{m\geq 0} C_{m}^{2}(2m)! \sup _{n} \frac{1}{n^{2}}\sum_{i,j\geq 1, i\not= j}^{n-1} R_{H_{1}}\left( \frac{i}{n}, \frac{j}{n} \right) ^{2m} \left( \frac{i}{n}\right) ^{-2H_{1}(m+\frac{1}{2})}\left( \frac{j}{n}\right) ^{-2H_{1}(m+\frac{1}{2})}\\
&=&\sum_{m\geq 0} C_{m}^{2}(2m)! \sup _{n}
\frac{1}{n^{2}}\sum_{i,j\geq 1, i\not= j}^{n-1} R_{H_{1}}\left( 1,
\left(\frac{j}{i}\right) \right) ^{2m} \left( \frac{j}{i} \right)
^{-2H_{1}m} \left( \frac{i}{n}\frac{j}{n}\right) ^{-H_{1}}.
\end{eqnarray*}
Let's focus on the case where $H_{1} < \frac{1}{2}$ first. Let
$Q_{H_{1}}\left(z\right)$ be the function defined by
\begin{eqnarray*}
Q_{H_{1}}\left(z\right) = \left\{ \begin{array}{ll}
\frac{R_{H_{1}}\left(1,z\right)}{z^{H_{1}}} \mbox{\   \   \ if\  \  \  } z \in \left(0,1\right]\\
0 \mbox{\   \   \ if\  \  \  } z = 0.
\end{array} \right.
\end{eqnarray*}
For $H_{1} < \frac{1}{2}$, we have
\begin{eqnarray*}
Q_{H_{1}}\left(z\right) \leq z^{H_{1}}.
\end{eqnarray*}
Indeed, the function $f(z) = 1 - z^{2H_{1}} - (1-z)^{2H_{1}}$ is
negative on $\left[0,1\right]$, increasing on
$\left[\frac{1}{2},1\right]$, decreasing on
$\left[0,\frac{1}{2}\right]$ and $f(1)=f(0)=0$. It follows that
\begin{eqnarray*}
\sum_{m\geq 0} C_{m}^{2}(2m)! \sup_{n,\varepsilon} \left|
A_{m}(\varepsilon ,n) \right| & \leq & 2\sum_{m\geq 0}
C_{m}^{2}(2m)! \sup _{n} \frac{1}{n^{2}}\sum_{i,j=0;i > j}^{n-1}
\left( \frac{j}{i} \right) ^{2H_{1}m} \left(
\frac{i}{n}\frac{j}{n}\right) ^{-H_{1}}
\\
& \leq & 2\sum_{m\geq 0} C_{m}^{2}(2m)! \sup _{n}
\frac{1}{n^{2}}\sum_{i,j=0;i > j}^{n-1} \left(\frac{j}{n} \right)
^{H_{1}(2m-1)} \left( \frac{i}{n}\right) ^{-H_{1}(2m+1)}
\\
& \leq & 2\sum_{m\geq 0} C_{m}^{2}(2m)! \sup _{n}
n^{2H_{1}-2}\sum_{i=0}^{n-1}i^{-H_{1}(2m+1)}\sum_{j=1}^{i-1}\int_{j}^{j+1}j^{H_{1}(2m-1)}dx
\\
& \leq & 2\sum_{m\geq 0} C_{m}^{2}(2m)! \sup _{n}
n^{2H_{1}-2}\sum_{i=0}^{n-1}i^{-H_{1}(2m+1)}\int_{0}^{i}x^{H_{1}(2m-1)}dx
\\
& \leq & 2\sum_{m\geq 0} C_{m}^{2}(2m)! \sup _{n}
\frac{n^{2H_{1}-2}}{2H_{1}m - H_{1} + 1}\sum_{i=0}^{n-1}i^{1-2H_{1}}
\\
& \leq & 2\sum_{m\geq 0} C_{m}^{2}(2m)! \sup _{n}
\frac{n^{-1}}{2H_{1}m - H_{1} + 1}\sum_{i=0}^{n-1}1
\\
& \leq & 2\sum_{m\geq 0} C_{m}^{2}(2m)! \sup _{n} \frac{1}{2H_{1}m -
H_{1} + 1}\leq  2\sum_{m\geq 0} \frac{C_{m}^{2}(2m)!}{2H_{1}m - H_{1} + 1}.
\end{eqnarray*}
Given that, by using Stirling's formula, the coefficient $C_{m}^{2}(2m)!$ behaves
as $\frac{1}{\sqrt{m}}$, we obtain that the above sum is finite.
Thus, we obtain the convergence of $a_{n}^{(3)}$ to zero in
$L^{2}(\Omega)$ for $H_{1}<\frac{1}{2}$.
\\\\
\noindent Let us now treat the case $H_{1}>\frac{1}{2}$.  We know (see \cite{Edd}, Lemma 1) that the function  $Q_{H}$ is increasing on $[0,1]$. Since $\frac{j}{i}\leq \frac{i-1}{i}= 1-\frac{1}{i}$ it holds that $Q_{H}(\frac{j}{i})\leq Q_{H}(1-\frac{1}{i})$. Then
\begin{eqnarray*}
\sum_{m\geq 0} C_{m}^{2}(2m)! \sup_{n,\varepsilon} \left|
A_{m}(\varepsilon ,n) \right| & \leq & 2\sum_{m\geq 0}
C_{m}^{2}(2m)! \sup _{n} \frac{1}{n^{2}}\sum_{i=1}^{n-1}Q_{H}\left( 1-\frac{1}{i}\right) \sum_{j=1}^{n-1}  \left(
\frac{i}{n}\frac{j}{n}\right) ^{-H_{1}}
\\
&=& 2\sum_{m\geq 0}
C_{m}^{2}(2m)! \sup _{n} \frac{1}{n}\sum_{i=1}^{n-1}Q_{H}\left( 1-\frac{1}{i}\right)\left(
\frac{i}{n}\right) ^{-H_{1}}\sum_{j=1}^{i-1} \int_{\frac{j-1}{n}}^{\frac{1}{n}}x^{-H_{1}}dx\\
&\leq & c_{H}\sum_{m\geq 0}
C_{m}^{2}(2m)! \sup _{n} \frac{1}{n}\sum_{i=1}^{n-1}Q_{H}\left( 1-\frac{1}{i}\right)\left(
\frac{i}{n}\right) ^{-H_{1}}\left(
\frac{i-1}{n}\right) ^{1-H_{1}}\\
&\sim & c_{H}\sum_{m\geq 0}
C_{m}^{2}(2m)! \sup _{n} \frac{1}{n}\sum_{i=1}^{n-1}Q_{H}\left( 1-\frac{1}{i}\right)\left(
\frac{i}{n}\right) ^{1-2H_{1}}.
\end{eqnarray*}
By adapting Lemma 2 in \cite{Edd} (by separating the sum over $i$ in a sum with $\frac{1}{i} \leq \delta$ and $\frac{1}{i}>\delta$ with $\delta $ suitably chosen),  we can prove that
\begin{equation*}
\frac{1}{n}\sum_{i,j=0}^{n-1} R_{H_{1}}\left( 1, \left(\frac{j}{i}\right) \right) ^{2m} \left( \frac{j}{i} \right) ^{-2H_{1}m} \left( \frac{i}{n}\frac{j}{n}\right) ^{-H}\leq c(H_{1}) m^{-\frac{1}{2H_{1}}}
\end{equation*}
with $c(H_{1})$ not depending on $m$ nor $n$. As a consequence
$$\sum_{m\geq 0} c_{m}^{2}(2m)!  \sup_{n,\eps} \left| A_{m}(\eps ,n) \right|\leq c(H_{1}) c_{m}^{2} (2m)! m^{-\frac{1}{2H_{1}}}.$$
The Stirling formula implies again that the above series is finite.\qed
\end{dem}

\begin{theorem}\label{t1}
Let $\langle S \rangle _{n}$ be given by (\ref{brac}). Then, as $n\to \infty$, we have the convergence in distribution
\begin{equation*}
n^{\alpha + H_{1} - 1}\langle  S \rangle _{n} \to \int_{\mathbb{R} } K^{2}(y)dy L^{H_{1}}(1,0)
\end{equation*}
where $L^{H_{1}}(1,0)$ is the local time of the fractional Brownian motion $B^{H_{1}}$.
\end{theorem}
{\bf Proof: } Using Proposition \ref{prop2} it suffices to check that $n^{\alpha +H_{1}} \int_{0}^{1} K^{2} (n^{\alpha +H_{1}} B^{H_{1}}_{s})ds$
converges to $\int_{\mathbb{R}}K^{2}(y)dy L^{H_{1}}(1,0) $. Using the occupation time formula (\ref{oc}), we obtain
\begin{equation*}
n^{\alpha +H_{1}} \int_{0}^{1} K^{2} (n^{\alpha +H_{1}} B^{H_{1}}_{s})ds= n^{\alpha + H_{1}} \int_{\mathbb{R}}K^{2}(n^{\alpha +H_{1}}x)L^{H_{1}}(1,x)dx = \int_{\mathbb{R}}K^{2}(y) L(1, yn^{-\alpha -H_{1}})dy
\end{equation*}
which converges as $n\to \infty$ to $\int_{\mathbb{R}}K^{2}(y)dy L^{H_{1}}(1,0) $ by using the continuity properties of the local time. \qed

\section{Limit distribution of $S_{n}$}
\noindent In this paragraph, we prove the limit in distribution of (\ref{sn}). Recall the notation (\ref{notationCovariancefBm}) and let's consider the Gaussian vector
$$X^{H_{2}} = (X_{1}^{H_{2}},...,X_{n}^{H_{2}}) = (B_{1}^{H_{2}} -
B_{0}^{H_{2}},...,B_{n}^{H_{2}}-B_{n-1}^{H_{2}}).$$ From this definition, it follows that
\begin{eqnarray*}
S_{n} = \sum_{i=0}^{n-1}K(n^{\alpha}B_{i}^{H_{1}}) (B_{i+1}^{H_{2}}
- B_{i}^{H_{2}}) =
\sum_{i=0}^{n-1}K(n^{\alpha}B_{i}^{H_{1}})X_{i+1}^{H_{2}}.
\end{eqnarray*}

\begin{theorem} Let $(S_{n})$ be given by (\ref{sn}) and assume that
\begin{equation} \label{cond2}
\alpha <1-H_{1}
\end{equation}
Then we have the convergence in law
\begin{equation*}
n^{\alpha +H_{1}-1} S_{n} \underset{n \rightarrow +\infty}{\longrightarrow} d_{1} W_{L^{H_{1}}(1,0)}
\end{equation*}
where $L^{H_{1}}(1,0)$ is the local time of $B^{H_{1}}$, $d_{1}:= \int _{\mathbb{R}}K^{2}(y)dy$ and  $W$ is a Brownian motion independent from $B^{H_{1}}$.
\end{theorem}
{\bf Proof: } We will study the characteristic function of $n^{\frac{\alpha}{2} + \frac{H_{1}}{2} -\frac{1}{2}}S_{n}$. In order to simplify the presentation, we will use the following notation. Let $i_{0}$ be the imaginary unit and $\lambda_{n}$ be
\begin{eqnarray*}
\lambda_{n} = \lambda n^{\frac{\alpha}{2} + \frac{H_{1}}{2} -\frac{1}{2}} \mbox{\  \  \  with \  \  }\lambda \in \mathbb{R}.
\end{eqnarray*}
Using the independence of the two fBms and computing the conditional expectation of $e^{i\lambda_{n} S_{n}}$ given $B^{H_{1}}$ we get
\begin{eqnarray*}
\mathbf{E}\left(e^{i_{0}\lambda_{n} S_{n}}\right) =
\mathbf{E}\left(e^{-\frac{1}{2}\sum_{i,j =
0}^{n-1}\lambda_{n}^{2}K\left(n^{\alpha}B_{i}^{H_{1}}\right)K\left(n^{\alpha}B_{j}^{H_{1}}\right)f_{H_{2}}(i,j)}\right)
\end{eqnarray*}
because if $X$ is a Gaussian vector with mean $\mu$ and covariance
matrix $\Sigma$, it's characteristic function is given by
$$\mathbf{E}\left(e^{i_{0}\left\langle t,X\right\rangle}\right) =
e^{i_{0}\mu^{\mbox{\tiny{T}}}t - \frac{1}{2}t^{\mbox{\tiny{T}}}\Sigma t}.$$ It follows that, with $f_{H_{2}}(i,j)$ given by (\ref{notationCovariancefBm}),
\begin{eqnarray*}
\mathbf{E}\left(e^{i_{0}\lambda_{n} S_{n}}\right) &=&
\mathbf{E}\left(e^{-\frac{\lambda_{n}^{2}}{2}\sum_{i =
0}^{n-1}K^{2}\left(n^{\alpha}B_{i}^{H_{1}}\right)}e^{-\frac{\lambda_{n}^{2}}{2}\sum_{i\neq
j =
0}^{n-1}K\left(n^{\alpha}B_{i}^{H_{1}}\right)K\left(n^{\alpha}B_{j}^{H_{1}}\right)f_{H_{2}}(i,j)}\right)
\\
&=& \mathbf{E}\left(e^{-\frac{\lambda_{n}^{2}}{2}\sum_{i =
0}^{n-1}K^{2}\left(n^{\alpha}B_{i}^{H_{1}}\right)}e^{-\lambda_{n}^{2}\sum_{i
= 0}^{n-1}\sum_{j =
0}^{i-1}K\left(n^{\alpha}B_{i}^{H_{1}}\right)K\left(n^{\alpha}B_{j}^{H_{1}}\right)f_{H_{2}}(i,j)}\right)
\\
&=& \mathbf{E}\left(e^{-\frac{\lambda_{n}^{2}}{2}\sum_{i =
0}^{n-1}K^{2}\left(n^{\alpha}B_{i}^{H_{1}}\right)}e^{-\lambda_{n}^{2}\sum_{i
= 0}^{n-1}\sum_{j =
0}^{i-1}K\left(n^{\alpha}B_{i}^{H_{1}}\right)K\left(n^{\alpha}B_{j}^{H_{1}}\right)H_{2}(2H_{2}-1)\int_{i}^{i+1}\int_{j}^{j+1}\left|s-u\right|^{2H_{2}-2}duds}\right)
\\
&=& \mathbf{E}\left(e^{-\frac{\lambda_{n}^{2}}{2}\sum_{i =
0}^{n-1}K^{2}\left(n^{\alpha}B_{i}^{H_{1}}\right)}e^{-\lambda_{n}^{2}H_{2}(2H_{2}-1)\int_{0}^{n}\int_{0}^{\left[s\right]}
K\left(n^{\alpha}B_{\left[s\right]}^{H_{1}}\right)K\left(n^{\alpha}B_{\left[u\right]}^{H_{1}}\right)\left|s-u\right|^{2H_{2}-2}duds}\right).
\end{eqnarray*}
Consider the process $(V_{n})_{n \geq 0}$ defined by $$V_{n} =
\int_{0}^{n}\int_{0}^{\left[s\right]}K\left(n^{\alpha}B_{\left[s\right]}^{H_{1}}\right)K\left(n^{\alpha}B_{\left[u\right]}^{H_{1}}\right)\left|s-u\right|^{2H_{2}-2}duds$$
and the function $\psi$ defined by $\psi(x) = e^{-\lambda_{n}^{2}H_{2}(2H_{2}-1)x}.$
Note that, since we excluded the diagonal, the integral $duds$ in the expression of $V_{n}$ makes sense even for $H_{2}<  \frac{1}{2}$. Note also  that $V_{n}$
is a bounded variation process (its quadratic variation is 0).
Furthermore, $$\psi'(x) =
-\lambda_{n}^{2}H_{2}(2H_{2}-1)e^{-\lambda_{n}^{2}H_{2}(2H_{2}-1)x}.
$$ Using the change of variables formula for bounded variation processes, it follows that
\begin{eqnarray*}
\psi(V_{n}) = 1 +\int_{0}^{n}\psi'(V_{s})dV_{s}
\end{eqnarray*}
i.e.,
\begin{eqnarray*}
e^{-\lambda_{n}^{2}H_{2}(2H_{2}-1)V_{n}} = 1
-\lambda_{n}^{2}H_{2}(2H_{2}-1)\int_{0}^{n}e^{-\lambda_{n}^{2}H_{2}(2H_{2}-1)V_{s}}dV_{s}.
\end{eqnarray*}
Therefore,
\begin{eqnarray*}
\mathbf{E}\left(e^{i_{0}\lambda_{n} S_{n}}\right) &=&
\mathbf{E}\left(e^{-\frac{\lambda_{n}^{2}}{2}\sum_{i =
0}^{n-1}K^{2}\left(n^{\alpha}B_{i}^{H_{1}}\right)}\left(1-\lambda_{n}^{2}H_{2}(2H_{2}-1)\int_{0}^{n}e^{-\lambda_{n}^{2}H_{2}(2H_{2}-1)V_{s}}dV_{s}\right)\right)
\\
&=& \mathbf{E}\left(e^{-\frac{\lambda_{n}^{2}}{2}\sum_{i =
0}^{n-1}K^{2}\left(n^{\alpha}B_{i}^{H_{1}}\right)}\right) \\
&& -
\mathbf{E}\left(\lambda_{n}^{2}H_{2}(2H_{2}-1)e^{-\frac{\lambda_{n}^{2}}{2}\sum_{i
=
0}^{n-1}K^{2}\left(n^{\alpha}B_{i}^{H_{1}}\right)}\int_{0}^{n}e^{-\lambda_{n}^{2}H_{2}(2H_{2}-1)V_{s}}dV_{s}\right)
\\
&:= & \mathbf{E}\left(T_{1}\right) - \mathbf{E}\left(T_{2}\right).
\end{eqnarray*}
We will now focus on the term $\mathbf{E}\left(T_{2}\right)$ and show that $$T_{2} \overset{L^{1}}{\longrightarrow} 0.$$
From
\begin{eqnarray*}
dV_{s} =
\left(\int_{0}^{\left[s\right]}K\left(n^{\alpha}B_{\left[s\right]}^{H_{1}}\right)K\left(n^{\alpha}B_{\left[u\right]}^{H_{1}}\right)\left|s-u\right|^{2H_{2}-2}du\right)ds
\end{eqnarray*}
we get
\begin{eqnarray*}
\mathbf{E}\left(T_{2}\right) &=&
\mathbf{E}\left(\lambda_{n}^{2}H_{2}(2H_{2}-1)e^{-\frac{\lambda_{n}^{2}}{2}\sum_{i
= 0}^{n-1}K^{2}\left(n^{\alpha}B_{i}^{H_{1}}\right)}\right.
\\
&& \left.\times
\int_{0}^{n}e^{-\lambda_{n}^{2}H_{2}(2H_{2}-1)V_{s}}\int_{0}^{\left[s\right]}K\left(n^{\alpha}B_{\left[s\right]}^{H_{1}}\right)K\left(n^{\alpha}B_{\left[u\right]}^{H_{1}}\right)\left|s-u\right|^{2H_{2}-2}duds\right)
\\
&=&
\mathbf{E}\left(\lambda_{n}^{2}H_{2}(2H_{2}-1)\int_{0}^{n}e^{-\frac{\lambda_{n}^{2}}{2}\int_{0}^{s}K^{2}\left(n^{\alpha}B_{\left[u\right]}^{H_{1}}\right)du}e^{-\frac{\lambda_{n}^{2}}{2}\int_{s}^{n}K^{2}\left(n^{\alpha}B_{\left[u\right]}^{H_{1}}\right)du}\right.
\\
&& \left.\times
e^{-\lambda_{n}^{2}H_{2}(2H_{2}-1)V_{s}}\int_{0}^{\left[s\right]}K\left(n^{\alpha}B_{\left[s\right]}^{H_{1}}\right)K\left(n^{\alpha}B_{\left[u\right]}^{H_{1}}\right)\left|s-u\right|^{2H_{2}-2}duds\right).
\end{eqnarray*}
Recall that the following holds
\begin{eqnarray}
\label{espCond}
\mathbf{E}\left(e^{i_{0}\lambda_{n} S_{s}}|B_{s}^{H_{1}}\right) =
\mathbf{E}\left(e^{-\frac{\lambda_{n}^{2}}{2}\int_{0}^{s}K^{2}\left(n^{\alpha}B_{\left[u\right]}^{H_{1}}\right)du}e^{-\lambda_{n}^{2}H_{2}(2H_{2}-1)V_{s}}|B_{s}^{H_{1}}\right).
\end{eqnarray}
This can be  seen for $s$ integer as at the beginning of this proof and also (\ref{espCond}) can easily be checked for any $s>0$.
We will use this property to compute the following upper bound for $\mathbf{E}\left(\left|T_{2}\right|\right)$
\begin{eqnarray*}
\mathbf{E}\left(\left|T_{2}\right|\right) & \leq &
\mathbf{E}\left(\lambda_{n}^{2}\int_{0}^{n}e^{-\frac{\lambda_{n}^{2}}{2}\int_{0}^{s}K^{2}\left(n^{\alpha}B_{\left[u\right]}^{H_{1}}\right)du}e^{-\lambda_{n}^{2}H_{2}(2H_{2}-1)V_{s}}\underbrace{\left|e^{-\frac{\lambda_{n}^{2}}{2}\int_{s}^{n}K^{2}\left(n^{\alpha}B_{\left[u\right]}^{H_{1}}\right)du}\right|}_{\leq
1}\right.
\\
&& \left. \times
\int_{0}^{\left[s\right]}K\left(n^{\alpha}B_{\left[s\right]}^{H_{1}}\right)K\left(n^{\alpha}B_{\left[u\right]}^{H_{1}}\right)H_{2}\left|2H_{2}-1\right|\left|s-u\right|^{2H_{2}-2}duds\right)
\\
& \leq &
\mathbf{E}\left(\lambda_{n}^{2}\int_{0}^{n}\mathbf{E}\left(e^{-\frac{\lambda_{n}^{2}}{2}\int_{0}^{s}K^{2}\left(n^{\alpha}B_{\left[u\right]}^{H_{1}}\right)du}e^{-\lambda_{n}^{2}H_{2}(2H_{2}-1)V_{s}}|B_{s}^{H_{1}}\right)\right.
\\
&& \left. \times
\int_{0}^{\left[s\right]}K\left(n^{\alpha}B_{\left[s\right]}^{H_{1}}\right)K\left(n^{\alpha}B_{\left[u\right]}^{H_{1}}\right)H_{2}\left|2H_{2}-1\right|\left|s-u\right|^{2H_{2}-2}duds\right).
\end{eqnarray*}
This is true because all the terms of the double integral are measurable with respect to the filtration generated by $(B_{u}^{H_{1}},u\leq s)$. At this point, we use (\ref{espCond}) to write
\begin{eqnarray*}
\mathbf{E}\left(\left|T_{2}\right|\right) & \leq &
\mathbf{E}\left(\lambda_{n}^{2}\int_{0}^{n}\mathbf{E}\left(e^{i_{0}\lambda_{n} S_{s}}|B_{s}^{H_{1}}\right)
\int_{0}^{\left[s\right]}K\left(n^{\alpha}B_{\left[s\right]}^{H_{1}}\right)K\left(n^{\alpha}B_{\left[u\right]}^{H_{1}}\right)H_{2}\left|2H_{2}-1\right|\left|s-u\right|^{2H_{2}-2}duds\right)
\\
& \leq &
\mathbf{E}\left(\lambda_{n}^{2}\int_{0}^{n}\underbrace{\left|e^{i_{0}\lambda_{n} S_{s}}\right|}_{=1} \int_{0}^{\left[s\right]}K\left(n^{\alpha}B_{\left[s\right]}^{H_{1}}\right)K\left(n^{\alpha}B_{\left[u\right]}^{H_{1}}\right)H_{2}\left|2H_{2}-1\right|\left|s-u\right|^{2H_{2}-2}duds\right)
\\
& \leq &
\mathbf{E}\left(\lambda_{n}^{2}\int_{0}^{n}\int_{0}^{\left[s\right]}K\left(n^{\alpha}B_{\left[s\right]}^{H_{1}}\right)K\left(n^{\alpha}B_{\left[u\right]}^{H_{1}}\right)H_{2}\left|2H_{2}-1\right|\left|s-u\right|^{2H_{2}-2}duds\right)
\\
& \leq &
\mathbf{E}\left(\lambda_{n}^{2}\sum_{i=0}^{n-1}\sum_{j=0}^{i-1}K\left(n^{\alpha}B_{i}^{H_{1}}\right)K\left(n^{\alpha}B_{j}^{H_{1}}\right)H_{2}\left|2H_{2}-1\right|\int_{i}^{i+1}\int_{j}^{j+1}\left|s-u\right|^{2H_{2}-2}duds\right).
\end{eqnarray*}
Assume that $H_{2} > \frac{1}{2}$, ergo
$\left|2H_{2}-1\right| > 0$ and $f_{H_{2}}(i,j) > 0$. Consequently,
\begin{eqnarray*}
\mathbf{E}\left(\left|T_{2}\right|\right) & \leq &
\mathbf{E}\left(\frac{\lambda_{n}^{2}}{2}\sum_{i=0}^{n-1}\sum_{j=0}^{i-1}K\left(n^{\alpha}B_{i}^{H_{1}}\right)K\left(n^{\alpha}B_{j}^{H_{1}}\right)f_{H_{2}}(i,j)\right)
\\
& \leq & \mathbf{E}\left(\frac{\lambda^{2}}{2}n^{\alpha + H_{1} -1}\sum_{i=0}^{n-1}\sum_{j=0}^{i-1}K\left(n^{\alpha}B_{i}^{H_{1}}\right)K\left(n^{\alpha}B_{j}^{H_{1}}\right)f_{H_{2}}(i,j)\right).
\end{eqnarray*}
The previous term is exactly the non-diagonal term of the $L^{2}$-norm of $n^{\frac{\alpha}{2} + \frac{H_{1}}{2} -\frac{1}{2}}S_{n}$ and we know that under condition (\ref{cond}), it converges to zero when $n \to +\infty$. Finally we have
$$\mathbf{E}\left(\left|T_{2}\right|\right) \underset{n \rightarrow
+\infty}{\longrightarrow} 0.$$ Assume now that $H_{2} < \frac{1}{2}$. It follows that $\left|2H_{2}-1\right| < 0$ and $f_{H_{2}}(i,j) < 0$, which gives us
\begin{eqnarray*}
\mathbf{E}\left(\left|T_{2}\right|\right) & \leq &
\mathbf{E}\left(-\frac{\lambda_{n}^{2}}{2}\sum_{i=0}^{n-1}\sum_{j=0}^{i-1}K\left(n^{\alpha}B_{i}^{H_{1}}\right)K\left(n^{\alpha}B_{j}^{H_{1}}\right)f_{H_{2}}(i,j)\right)
\\
& \leq & \mathbf{E}\left(-\frac{\lambda^{2}}{2}n^{\alpha + H_{1} -1}\sum_{i=0}^{n-1}\sum_{j=0}^{i-1}K\left(n^{\alpha}B_{i}^{H_{1}}\right)K\left(n^{\alpha}B_{j}^{H_{1}}\right)f_{H_{2}}(i,j)\right).
\end{eqnarray*}
As in the previous case, this term is again exactly the non-diagonal term of the $L^{2}$-norm of $n^{\frac{\alpha}{2} + \frac{H_{1}}{2} -\frac{1}{2}}S_{n}$ and for the same reasons, we get the following result again (which is now valid for any $H_{2} \in \left(0,1\right)$) $$\mathbf{E}\left(\left|T_{2}\right|\right) \underset{n \rightarrow
+\infty}{\longrightarrow} 0.$$
Concerning the term $T_{1}$, we note that
$$\mathbf{E}\left(T_{1}\right) =\mathbf{E}\left(e^{-\frac{\lambda ^{2} }{2}\langle S\rangle _{n}}\right)$$
and the result follows from Theorem \ref{t1}. \qed

\vskip0.3cm

\begin{remark}
The following comments deal with the conditions (\ref{cond}) and (\ref{cond2}). Condition (\ref{cond2}) is a natural extension of the condition  $\alpha <\frac{1}{2}$ in e.g. \cite{WP1}, \cite{WP2} which means that the bandwidth parameter satisfies $nh^{2}_{n} =nn^{-2\alpha } \to \infty$ as $n\to \infty$.
From (\ref{cond}) and (\ref{cond2}), this is the constraint we find for $\alpha$ (considering $\alpha$ is our degree of freedom)
$$
\left\{
\begin{array}{rl}
0<H_{1}<1 \\
0<H_{2}<1 \\
\alpha > 4H_{2} - H_{1} - 2 \\
\alpha < 1 - H_{1}
\end{array}
\right. \Leftrightarrow \left\{
\begin{array}{rl}
0<H_{1}<1 \\
0<H_{2}<1 \\
4H_{2} - H_{1} - 2 < \alpha < 1 - H_{1}.
\end{array}
\right.
$$ As an example, consider the case where $H_{1} = H_{2} = H$. Those constraints become
$$
\left\{
\begin{array}{rl}
0<H<1 \\
\left(3H - 2\right)^{+} < \alpha < 1 - H.
\end{array}
\right.
$$ For this system to have a solution, we need to verify that $$3H - 2 < 1 - H \Leftrightarrow H < \frac{3}{4}.$$ As a result, our constraints become
$$
\left\{
\begin{array}{rl}
0<H<\frac{3}{4} \\
\left(3H - 2\right)^{+} < \alpha < 1 - H.
\end{array}
\right.
$$
We could also consider the case where $\alpha$ has a fixed value and where the constraints would be on $H_{1}$ and $H_{2}$.
\end{remark}

\vskip0.5cm

\section{The stable convergence}
In this section we will study the convergence of the vector $(S_{n}, (G_{t}) _{t\geq 0} )$ where $(G_{t})_{t\geq 0}$ is a stochastic process independent from $B^{H_{1}}$ and satisfies some additional conditions. In this case, since the process $(G_{t})_{t\geq 0}$ is not necessarily a Gaussian process and since no information is available on the correlation between $B^{H_{2}}$ and $G_{t}$, the characteristic function of the vector  $(S_{n}, (G_{t}) _{t\geq 0} )$ cannot be computed directly. To compute it, we will use the tools of the stochastic calculus with respect to the fractional Brownian motion. The basic observation is that $S_{n}$ can be expressed as a stochastic integral with respect to $B^{H_{2}}$. Indeed,
\begin{eqnarray}
S_{n} & = & \sum_{i=0}^{n-1}K(n^{\alpha}B_{i}^{H_{1}}) (B_{i+1}^{H_{2}} - B_{i}^{H_{2}})
 =  \sum_{i=0}^{n-1}K(n^{\alpha}B_{i}^{H_{1}}) \delta^{H_{2}}(\textbf{1}_{\left[i,i+1\right]}(\cdot))\\
 & = & \sum_{i=0}^{n-1}\delta^{H_{2}}(K(n^{\alpha}B_{i}^{H_{1}}) \textbf{1}_{\left[i,i+1\right]}(\cdot)) + \left\langle \underbrace{D^{H_{2}}K(n^{\alpha}B_{i}^{H_{1}})}_{= 0\mbox{\  \footnotesize{from}\  }B^{H_{1}} \bot B_{t}^{H_{2}} },\textbf{1}_{\left[i,i+1\right]}(\cdot)\right\rangle_{\mathcal{H}_{H_{2}}} \nonumber
\\
& = & \sum_{i=0}^{n-1}\int_{i}^{i+1}K(n^{\alpha}B_{i}^{H_{1}})dB_{s}^{H_{2}} =
 \sum_{i=0}^{n-1}\int_{i}^{i+1}K(n^{\alpha}B_{\left[s\right]}^{H_{1}})dB_{s}^{H_{2}}  =  \int_{0}^{n}K(n^{\alpha}B_{\left[s\right]}^{H_{1}})dB_{s}^{H_{2}}. \nonumber
\end{eqnarray}
 We will also use the ``bracket'' of $S_{n}$. This quantity equals
\begin{eqnarray}
\sum_{i=0}^{n-1}K^{2}(n^{\alpha}B_{i}^{H_{1}}) & = & \sum_{i=0}^{n-1}\int_{i}^{i+1}K^{2}(n^{\alpha}B_{i}^{H_{1}})ds \nonumber
\\
& = & \sum_{i=0}^{n-1}\int_{i}^{i+1}K^{2}(n^{\alpha}B_{\left[s\right]}^{H_{1}})ds = \int_{0}^{n}K^{2}(n^{\alpha}B_{\left[s\right]}^{H_{1}})ds. \nonumber
\end{eqnarray}
Before going any further, we will describe the elements of the stochastic calculus with respect to fractional Brownian motion that we will be using in the sequel. We will start by introducing some notations and definitions. Let $\phi$ be the function defined by $$\phi(s,t) = H(2H-1)\left|s-t\right|^{2H-2}, \mbox{\  \  }s,t \in \mathbb{R}.$$ Let $D$ (introduced in section 2) be the Malliavin derivative operator with respect to the fractional Brownian motion with Hurst parameter $H$. Based on this operator, let $D^{\phi}$ be another derivative operator (called the $\phi$-derivative operator) defined by $$D_{t}^{\phi}F = \int_{\mathbb{R}}\phi(t,v)D_{v}F dv$$ for any $F$ in the domain of $D$. For more details about this operator, see \cite{Bia}. Let $\mathcal{L}_{\phi}(0,T)$ be the family of stochastic processes $F$ on $\left[0,T\right]$ with the following properties: $F \in \mathcal{L}_{\phi}(0,T)$ if and only if $\mathbf{E}\left[\left\|F\right\|_{\mathcal{H}}^{2}\right] < \infty$, $F$ is $\phi$-differentiable, the trace of $D_{s}^{\phi}F_{t},0\leq s,t\leq T$, exists, and $\mathbf{E}\left[\int_{0}^{T}\int_{0}^{T}\left|D_{s}^{\phi}F_{t}\right|^{2}dsdt\right]<\infty$ and for each sequence of partitions $(\pi_{n},n \in \mathbb{N})$ such that $\left|\pi_{n}\right|\rightarrow 0$ as $n \rightarrow +\infty$, $$\sum_{i,j=0}^{n-1}\mathbf{E}\left[\int_{t_{i}^{(n)}}^{t_{i+1}^{(n)}}\int_{t_{j}^{(n)}}^{t_{j+1}^{(n)}}\left|D_{s}^{\phi}F_{t_{i}^{(n)}}^{\pi}D_{t}^{\phi}F_{t_{j}^{(n)}}^{\pi} - D_{s}^{\phi}F_{t}D_{t}^{\phi}F_{s}\right|dsdt\right]$$ and $$\mathbf{E}\left[\left\|F^{\pi} - F\right\|_{\mathcal{H}}^{2}\right]$$ tend to 0 as $n \rightarrow +\infty$, where $\pi_{n}:<0=t_{0}^{(n)}<t_{1}^{(n)}<...<t_{n-1}^{(n)}<t_{n}^{(n)} = T$.\\\\
\noindent In our particular situation, we are dealing with processes of the form $\int_{0}^{t}F_{u}dB_{u}^{H} + \int_{0}^{t}G_{u}du$, (where $B^{H}$ is a fractional Brownian motion with Hurst parameter $H$), for which the following It\^o formula holds in the case $H>\frac{1}{2}$.
\begin{theorem}
Let $\eta_{t} = \int_{0}^{t}F_{u}dB_{u}^{H} + \int_{0}^{t}G_{u}du$, for $t \in \left[0,T\right]$ with $\mathbf{E}\left[\underset{0 \leq s \leq T}{\mbox{sup}}\left|G_{s}\right|\right] < \infty$ and let $(F_{u}, 0 \leq u \leq T)$ be a stochastic process in $\mathcal{L}_{\phi}(0,T)$. Assume that there is a $\beta > 1 - H$ such that
\begin{eqnarray}
\mathbf{E}\left[\left|F_{u} - F_{v}\right|^{2}\right] \leq C\left|u-v\right|^{2\beta}
\end{eqnarray}
where $\left|u-v\right| \leq \zeta$ for some $\zeta > 0$ and
\begin{eqnarray}
\underset{0 \leq u, v\leq t,\left|u-v\right| \rightarrow 0 }{\mbox{lim}}\mathbf{E}\left[\left|D_{u}^{\phi}(F_{u} - F_{v})\right|^{2}\right] = 0.
\end{eqnarray}
Let $f:\mathbb{R}_{+} \times \mathbb{R} \rightarrow \mathbb{R}$ be a function having the first continous derivative in its first variable and the second continous derivative in its second variable. Assume that these derivatives are bounded. Moreover, it is assumed that $\mathbf{E}\left[\int_{0}^{T}\left|F_{s}D_{s}^{\phi}\eta_{s}\right|ds\right] < \infty$ and $(\frac{\partial f(s,\eta_{s})}{\partial x}F_{s},s \in \left[0,T\right]) \in \mathcal{L}_{\phi}(0,T)$. Then for $t \in \left[0,T\right]$,
\begin{eqnarray}
f(t,\eta_{t}) & = & f(0,0) + \int_{0}^{t}\frac{\partial f}{ \partial s}(s,\eta_{s})ds + \int_{0}^{t}\frac{\partial f}{ \partial x}(s,\eta_{s})G_{s}ds \nonumber
\\
& & + \int_{0}^{t}\frac{\partial f}{ \partial x}(s,\eta_{s})F_{s}dB_{s}^{H} + \int_{0}^{t}\frac{\partial^{2} f}{ \partial x^{2}}(s,\eta_{s})F_{s}D_{s}^{\phi}\eta_{s}ds.
\end{eqnarray}
\end{theorem}
We also have the following technical lemma (\cite{Bia} p.71) that will be particularly useful for our future computations.
\begin{lemma}
\label{tekLem}
Let $(F_{t},t \in \left[0,T\right])$  be a stochastic process in $\mathcal{L}_{\phi}(0,T)$ and $$\underset{0\leq s \leq T }{\mbox{sup}}\mathbf{E}\left[\left|D_{s}^{\phi}F_{s}\right|^{2}\right] < \infty$$ and let $\eta_{t} = \int_{0}^{t}F_{u}dB_{u}^{H}$ for $t \in \left[0,T\right]$. Then, for $s,t \in \left[0,T\right]$,
\begin{eqnarray}
D_{s}^{\phi}\eta_{t} = \int_{0}^{t}D_{s}^{\phi}F_{u}dB_{u}^{H} + \int_{0}^{t}F_{u}\phi(s,u)du.
\end{eqnarray}
\end{lemma}
It is now possible to state the main result of this section
\begin{theorem}
Assume that (\ref{cond}) and (\ref{cond2}) holds. Let $(G_{t})_{t\geq 0} $ be a stochastic process independent from $B^{H_{1}}$ and adapted to the filtration generated by $B^{H_{2}}$ such that for every $t\geq 0$ the random variable $G_{t}$ belongs to
$\mathbb{D}^{1,2}$ and $\Vert D_{s}G_{t}\vert \leq C$ for any $s,t$ and $\omega$. Then the vector $(S_{n}, (G_{t})_{t\geq 0})$ converges in the sense of finite dimensional distributions to the vector $(cW_{L^{H_{1}}(1,0)}, (G_{t})_{t\geq 0}$, where $c$ is a positive constant.
\end{theorem}
{\bf Proof: } In order to simplify the presentation, the following notations will be used. We will denote by $\lambda_{n}$ (like we did in a previous proof) the quantity $$\lambda_{n} = \lambda n^{\frac{\alpha}{2} + \frac{H_{1}}{2} - \frac{1}{2}}$$ where $\lambda \in \mathbb{R}$. The following notation will also be used:
\begin{equation*}
e(\lambda ,n)= e^{-\frac{\lambda^{2}}{2}n^{\alpha + H_{1} - 1}\int_{0}^{n}K^{2}(n^{\alpha}B_{\left[u\right]}^{H_{1}})du}=e^{-\frac{\lambda^{2}}{2}n^{\alpha + H_{1} - 1}\sum_{i=0}^{n-1}K^{2}(n^{\alpha } B^{H_{1}}_{i})}.
\end{equation*}
Let $(F_{t},t \geq 0)$ and $(G_{t},t \geq 0)$ be two stochastic processes defined by
$$
\left\{
\begin{array}{ll}
F_{u} = K(n^{\alpha}B_{\left[u\right]}^{H_{1}}) \\
G_{u} = -i_{0}\frac{\lambda_{n}}{2}K^{2}(n^{\alpha}B_{\left[u\right]}^{H_{1}})
\end{array}
\right.
$$
and let $(\eta_{t}^{(\lambda_{n})},t \geq 0)$ be the stochastic process defined by
\begin{eqnarray*}
\eta_{t}^{(\lambda_{n})} = \int_{0}^{t}F_{u}dB_{u}^{H} + \int_{0}^{t}G_{u}du =  \int_{0}^{t}K(n^{\alpha}B_{\left[u\right]}^{H_{1}})dB_{u}^{H_{2}} -i_{0}\frac{\lambda_{n}}{2} \int_{0}^{t}K^{2}(n^{\alpha}B_{\left[u\right]}^{H_{1}})du.
\end{eqnarray*}
Consider the function $f:\mathbb{C} \rightarrow \mathbb{C}, f(x) = e^{i_{0}\lambda_{n} x}$. We can apply the It\^o formula to $f(\eta_{t}^{(\lambda_{n})})$ in order to obtain
\begin{eqnarray}
\label{itoForm}
e^{i_{0}\lambda_{n} \eta_{t}^{(\lambda_{n})}} & = & 1 + \frac{\lambda_{n}^{2}}{2}\int_{0}^{t} e^{i_{0} \lambda_{n} \eta_{s}^{(\lambda_{n})}}K^{2}(n^{\alpha}B_{\left[s\right]}^{H_{1}})ds \nonumber
\\
&& + i_{0}\lambda_{n} \int_{0}^{t} e^{i_{0} \lambda_{n} \eta_{s}^{(\lambda_{n})}}K(n^{\alpha}B_{\left[s\right]}^{H_{1}})dB_{s}^{H_{2}} \nonumber
\\
&& -\lambda_{n}^{2}\int_{0}^{t} e^{i_{0} \lambda_{n} \eta_{s}^{(\lambda_{n})}}K(n^{\alpha}B_{\left[s\right]}^{H_{1}})D_{s}^{\phi,H_{2}}\eta_{s}^{(\lambda_{n})}ds
\end{eqnarray}
where $D^{\phi,H_{2}}$ is the operator $D^{\phi}$ introduced above with respect to the fractional Brownian motion $B^{H_{2}}$.  We use Lemma \ref{tekLem} to compute $D_{s}^{\phi,H_{2}}\eta_{s}^{(\lambda_{n})}$. We get
\begin{eqnarray}
D_{s}^{\phi,H_{2}}\eta_{s}^{(\lambda_{n})} & = & D_{s}^{\phi,H_{2}}\int_{0}^{s}K(n^{\alpha}B_{\left[u\right]}^{H_{1}})dB_{u}^{H_{2}} -i_{0}\frac{\lambda_{n}}{2} \underbrace{D_{s}^{\phi,H_{2}}\int_{0}^{s}K^{2}(n^{\alpha}B_{\left[u\right]}^{H_{1}})du}_{= 0\mbox{\  \footnotesize{from}\  }B^{H_{1}} \bot B_{t}^{H_{2}} } \nonumber
\\
& = & \underbrace{\int_{0}^{s}D_{s}^{\phi,H_{2}}K(n^{\alpha}B_{\left[u\right]}^{H_{1}})dB_{u}^{H_{2}}}_{= 0\mbox{\  \footnotesize{from}\  }B^{H_{1}} \bot B_{t}^{H_{2}} } + H_{2}(2H_{2} - 1) \int_{0}^{s}K(n^{\alpha}B_{\left[u\right]}^{H_{1}})\left|s-u\right|^{2H_{2} - 2}du \nonumber
\\
& = & H_{2}(2H_{2} - 1) \int_{0}^{s}K(n^{\alpha}B_{\left[u\right]}^{H_{1}})\left|s-u\right|^{2H_{2} - 2}du. \nonumber
\end{eqnarray}
By substituting in (\ref{itoForm}), we obtain
\begin{eqnarray}
\label{itoExpanded}
e^{i_{0}\lambda_{n} \eta_{t}^{(\lambda_{n})}} & = & 1 + \frac{\lambda_{n}^{2}}{2}\int_{0}^{t} e^{i_{0} \lambda_{n} \eta_{s}^{(\lambda_{n})}}K^{2}(n^{\alpha}B_{\left[s\right]}^{H_{1}})ds \nonumber
\\
&& + i_{0}\lambda_{n} \int_{0}^{t} e^{i_{0} \lambda_{n} \eta_{s}^{(\lambda_{n})}}K(n^{\alpha}B_{\left[s\right]}^{H_{1}})dB_{s}^{H_{2}} \nonumber
\\
&& -\lambda_{n}^{2}H_{2}(2H_{2} - 1)\int_{0}^{t} e^{i_{0} \lambda_{n} \eta_{s}^{(\lambda_{n})}}K(n^{\alpha}B_{\left[s\right]}^{H_{1}}) \nonumber
\\
&& \times \int_{0}^{s}K(n^{\alpha}B_{\left[u\right]}^{H_{1}})\left|s-u\right|^{2H_{2} - 2}duds.
\end{eqnarray}
By multiplying both sides of (\ref{itoExpanded}) by $e^{-\frac{\lambda_{n}^{2}}{2}\int_{0}^{n}K^{2}(n^{\alpha}B_{\left[u\right]}^{H_{1}})du} = e(\lambda,n)$, we obtain
\begin{eqnarray}
&& e^{i_{0}\lambda_{n}\int_{0}^{n}K(n^{\alpha}B_{\left[u\right]}^{H_{1}})dB_{u}^{H_{2}} }
 = e^{-\frac{\lambda_{n}^{2}}{2}\int_{0}^{n}K^{2}(n^{\alpha}B_{\left[u\right]}^{H_{1}})du} \nonumber
\\
&& + \frac{\lambda_{n}^{2}}{2}\int_{0}^{n} e^{i_{0} \lambda_{n} \eta_{s}^{(\lambda_{n})}}K^{2}(n^{\alpha}B_{\left[s\right]}^{H_{1}})ds \cdot e(\lambda,n) \nonumber
\\
&& + i_{0}\lambda_{n} \int_{0}^{n} e^{i_{0} \lambda_{n}
\eta_{s}^{(\lambda_{n})}}K(n^{\alpha}B_{\left[s\right]}^{H_{1}})dB_{s}^{H_{2}}\cdot e(\lambda,n)\nonumber
\\
&& -\lambda_{n}^{2} H_{2}(2H_{2} - 1)\int_{0}^{n} e^{i_{0} \lambda_{n}  \eta_{s}^{(\lambda_{n})}}K(n^{\alpha}B_{\left[s\right]}^{H_{1}})\nonumber
\\
&& \times \int_{0}^{s}K(n^{\alpha}B_{\left[u\right]}^{H_{1}})\left|s-u\right|^{2H_{2} - 2}duds\cdot e(\lambda,n).\label{1000}
\end{eqnarray}
The sum of the last two terms in (\ref{1000}) can be written in a more suitable way by using sums instead of integrals. Together, these two last terms give us
\begin{eqnarray}
 & & \mathbf{E}\left(\lambda_{n}^{2}\int_{0}^{n}e^{i_{0}\lambda_{n}\eta_{s}^{(\lambda_{n})}}K\left(n^{\alpha}B_{\left[s\right]}^{H_{1}}\right)\right.\nonumber
\\
& & \left.\times \left[\frac{1}{2}K\left(n^{\alpha}B_{\left[s\right]}^{H_{1}}\right) - H_{2}(2H_{2} - 1)\int_{0}^{s}K\left(n^{\alpha}B_{\left[u\right]}^{H_{1}}\right)\left|s-u\right|^{2H_{2} - 2}du\right]ds \cdot e(\lambda , n)\right)\nonumber
\\
& = & \mathbf{E}\left(\lambda_{n}^{2}\int_{0}^{n}e^{i_{0}\lambda_{n}\eta_{s}^{(\lambda_{n})}}K\left(n^{\alpha}B_{\left[s\right]}^{H_{1}}\right)\right.\nonumber
\\
& & \left.\times \left[\frac{1}{2}K\left(n^{\alpha}B_{\left[s\right]}^{H_{1}}\right) - H_{2}(2H_{2} - 1)\left(\int_{0}^{i}K\left(n^{\alpha}B_{\left[u\right]}^{H_{1}}\right)\left|s-u\right|^{2H_{2} - 2}du\right.\right.\right.\nonumber
\\
& & \left.\left.\left.+ \int_{i}^{s}K\left(n^{\alpha}B_{i}^{H_{1}}\right)\left|s-u\right|^{2H_{2} - 2}du\right)\right]ds \cdot e(\lambda ,n)\right)\nonumber
\\
& = &-\mathbf{E}\left( \lambda_{n}^{2}\sum_{i = 0}^{n-1}\int_{i}^{i+1}e^{i_{0}\lambda_{n} \eta_{s}^{(\lambda_{n})}}K\left(n^{\alpha}B_{i}^{H_{1}}\right)\right.\nonumber
\\
& & \left.\times  H_{2}(2H_{2} - 1)\sum_{j = 0}^{i - 1}K\left(n^{\alpha}B_{j}^{H_{1}}\right)\int_{j}^{j+1}\left|s-u\right|^{2H_{2} - 2}duds \cdot e(\lambda ,n)\right)\nonumber
\\
& & + \mathbf{E}\left(\lambda_{n}^{2}\sum_{i = 0}^{n-1}\int_{i}^{i+1}e^{i_{0}\lambda_{n} \eta_{s}^{(\lambda_{n})}}K^{2}\left(n^{\alpha}B_{i}^{H_{1}}\right)\right. \nonumber
\\
& & \left.\times \left[\frac{1}{2} - H_{2}(2H_{2} - 1)
\int_{i}^{s}\left|s-u\right|^{2H_{2} - 2}du\right]ds\cdot e(\lambda
,n)\right). \label{1001}
\end{eqnarray}
Let us now fix $\beta _{1}, \ldots , \beta_{N}\in \mathbb{R}$ and $t_{1},\ldots , t_{N} \geq 0$. We need to show that $\mathbf{E}\left(
e^{i_{0}\lambda_{n}S_{n}}
e^{i_{0}(\beta_{1}G_{t_{1}}+\ldots + \beta_{N}G_{t_{N}})}\right)$ converges to $\mathbf{E}\left( e^{-\frac{\lambda
^{2}(L^{H_{1}}(1,0))^{2}}{2}}e^{i_{0}(\beta_{1}G_{t_{1}}+\ldots + \beta_{N}G_{t_{N}})}\right)$. We will use the notation
\begin{equation*}
g_{N}:=e^{i_{0}(\beta_{1}G_{t_{1}}+\ldots + \beta_{N}G_{t_{N}})}.
\end{equation*}
By combining relations (\ref{1000}) and (\ref{1001}), we can write
\begin{eqnarray}
\mathbf{E}\left( e^{i_{0}\lambda_{n}S_{n}} g_{N}\right)&=&\mathbf{E}(e(\lambda ,n)
g_{N})\nonumber
\\
&&+ \mathbf{E}\left(i_{0}\lambda _{n}\int_{0}^{n} e^{i_{0} \lambda_{n}
\eta_{s}^{(\lambda_{n})}}K(n^{\alpha}B_{\left[s\right]}^{H_{1}})dB_{s}^{H_{2}}\cdot
e(\lambda , n)g_{N}\right)\nonumber
\\
&&-\mathbf{E} \left(\lambda_{n} ^{2}\sum_{i = 0}^{n-1}\int_{i}^{i+1}e^{i_{0}\lambda
_{n}\eta_{s}^{(\lambda_{n})}}K\left(n^{\alpha}B_{i}^{H_{1}}\right)\right.
\nonumber
\\
& & \left.\times  H_{2}(2H_{2} - 1)\sum_{j = 0}^{i -
1}K\left(n^{\alpha}B_{j}^{H_{1}}\right)\int_{j}^{j+1}
\left|s-u\right|^{2H_{2} - 2}duds \times e(\lambda ,n)g_{N}\right)\nonumber
\\
& & + \mathbf{E}\left(\lambda _{n} ^{2}\sum_{i = 0}^{n-1}\int_{i}^{i+1}e^{i_{0}\lambda
_{n} \eta_{s}^{(\lambda_{n})}}K^{2}\left(n^{\alpha}B_{i}^{H_{1}}\right)\right.\nonumber
\\
& & \left.\times \left[\frac{1}{2} - H_{2}(2H_{2} - 1)
\int_{i}^{s}\left|s-u\right|^{2H_{2} - 2}du\right]ds\times e(\lambda
,n)g_{N}\right)\nonumber
\\
&: =& \mathbf{E}(e(\lambda ,n) g_{N})+T_{1}^{\star}+T_{2}^{\star}+
T_{3}^{\star}. \label{1002}
\end{eqnarray}
Let us begin by proving that the term $T_{2}^{\star}$ converges to zero as $n\to \infty$. Since
\begin{eqnarray}
\label{bound1}
\left| e^{i_{0}\lambda_{n}\eta_{u}^{(\lambda_{n})}}\right|e^{-\frac{\lambda_{n}^{2}}{2}\int_{0}^{n}K^{2}\left(n^{\alpha}B_{\left[u\right]}^{H_{1}}\right)du}
\leq 1
\end{eqnarray}
for every $s\leq n$ and since $\vert e^{i_{0}x}\vert =1$ for every $x$
real, $T_{2}^{\star}$ can be bounded as follows
\begin{eqnarray*}
T_{2}^{\star} &\leq &\mathbf{E}\left( \lambda_{n}^{2}c(H_{2}) \sum_{i =
0}^{n-1}\sum_{j=0}^{i-1}K\left(n^{\alpha}B_{i}^{H_{1}}\right)K\left(n^{\alpha}B_{j}^{H_{1}}\right)\int_{i}^{i+1}\int_{j}^{j+1}\vert s-u\vert ^{2H_{2}-2} duds\right)
\\
&\leq &\mathbf{E}\left( \lambda^{2} n^{\alpha + H_{1} - 1}c(H_{2}) \sum_{i =
0}^{n-1}\sum_{j=0}^{i-1}K\left(n^{\alpha}B_{i}^{H_{1}}\right)K\left(n^{\alpha}B_{j}^{H_{1}}\right)\int_{i}^{i+1}\int_{j}^{j+1}\vert s-u\vert ^{2H_{2}-2} duds\right)
\end{eqnarray*}
and this goes to zero as in the proof showing that the non-diagonal term goes to zero under the renormalization $n^{\alpha+
H_{1}-1}$. Let us now handle the term $T_{1}^{\star}$. By using the independence of $B^{H_{1}}$ and $B^{H_{2}}$ we can write \begin{equation*}
T_{1}^{\star}=\mathbf{E}\left(i_{0}\lambda_{n} \int_{0}^{n} e^{i_{0} \lambda_{n}\eta_{s}^{(\lambda_{n})}}K(n^{\alpha}B_{\left[s\right]}^{H_{1}})e(\lambda , n)
dB_{s}^{H_{2}}\cdot g_{N}\right).
\end{equation*}
The duality formula is used to obtain
\begin{eqnarray*}
T_{1}^{\star}&=&\mathbf{E}\left(i_{0}\lambda_{n}\langle \mathbf{1}_{\left[0,n\right]} e^{i_{0}\lambda_{n}\eta_{\cdot}^{(\lambda_{n})} }
K(n^{\alpha} B^{H_{1}}_{[\cdot ]})e(\lambda ,n), D^{H_{2}}g_{N} \rangle _{
{\cal{H}}_{H_{2}}}\right)\\
&=&\mathbf{E}\left(-\lambda_{n}\langle \mathbf{1}_{\left[0,n\right]} e^{i_{0}\lambda_{n}\eta_{\cdot}^{(\lambda_{n})} }
K(n^{\alpha} B^{H_{1}}_{[\cdot ]})e(\lambda ,n), g_{N} \sum_{k=1}^{N}\beta_{k} D^{H_{2}} G_{t_{k}}\rangle _{ {\cal{H}}_{H_{2}}}\right).
\end{eqnarray*}
Recall that the following formula holds (see \cite{N} for further details)
\begin{eqnarray*}
\left\langle \phi,\psi \right\rangle_{{\cal{H}}_{H_{2}}} = H_{2}(2H_{2}-1)\int_{0}^{T}\int_{0}^{T}\left|r-u\right|^{2H_{2}-2}\phi_{r}\psi_{u}dudr
\end{eqnarray*}
for any pair of functions in the Hilbert space ${\cal{H}}_{H_{2}}$. This formula is used to write $T_{1}^{\star}$ as
\begin{eqnarray*}
T_{1}^{\star}&=&\mathbf{E}\left(-\lambda_{n}\sum_{k=1}^{N}\beta_{k}H_{2}(2H_{2}-1)\int_{0}^{n}
\int_{0}^{t_{k}}e^{i_{0}\lambda_{n} \eta_{u}}K(n^{\alpha} B^{H_{1}}_{[u]} ) e(\lambda ,n)D_{v}G_{t_{k}}\vert u-v\vert
^{2H_{2}-2}dvdu\right)
\end{eqnarray*}
where the fact that $G_{t}$ is adapted to the filtration of $B^{H_{2}}$ is used. It suffices to show that for every fixed $t\geq 0$,
$$\lambda _{n}\mathbf{E}\left(\int_{0}^{n}
\int_{0}^{t}e^{i_{0}\lambda _{n} \eta _{u} }K(n^{\alpha} B^{H_{1}}_{[u]}
) e(\lambda ,n)D_{v}G_{t}\vert u-v\vert ^{2H_{2}-2}dvdu\right)$$ converges
to zero as $n\to \infty.$ Since  the derivative of $G_{t}$ is
bounded and using (\ref{bound1}) we find that  the above term is less than
\begin{eqnarray*}
&&\lambda _{n} c_{1} \mathbf{E}\left(\int_{0}^{n}\int_{0}^{t}K(n^{\alpha}
B^{H_{1}}_{[u]})\vert u-v\vert ^{2H_{2}-2}dvdu\right)
\\
&=& \lambda_{n}c_{1}\sum_{i=0}^{n-1} \mathbf{E}\left(K(n^{\alpha} B^{H_{1}}_{i})\right)
\int_{i}^{i+1}\int_{0}^{t}\vert u-v\vert ^{2H_{2}-2}dvdu
\\
&=& \lambda _{n} c_{1}c_{H_{2}}\sum_{i=0}^{n-1} \mathbf{E}\left(K(n^{\alpha}
B^{H_{1}}_{i})\right) \left( -\vert i+1-t\vert ^{2H_{2}} +\vert i-t\vert
^{2H_{2}} +\vert i+1\vert ^{2H_{2}}-i^{2H_{2}}\right)
\end{eqnarray*}
where $c_{1}$ is the constant upper bound of the derivative of $G_{t}$ and $c_{H_{2}}$ is a constant depending only on $H_{2}$. Since for every fixed $t>0$ the function $\left( -\vert i+1-t\vert ^{2H_{2}} +\vert i-t\vert ^{2H_{2}} +\vert i+1\vert
^{2H_{2}}-i^{2H_{2}}\right)=$ behaves, modulo a constant,  as $i^{2H_{2}-2}$ and since the order of the expectation of $K(n^{\alpha} B^{H_{1}}_{i})$ is the same as that of $n^{-\alpha} i^{-H_{1}}$ it is clear that $T_{1}$ converges to zero as $n\to \infty$.
\\\\
\noindent Finally, we will show that $T_{3}^{\star}$ converges to zero. Note that the term $T_{3}^{\star}$ can be expressed as follows
\begin{eqnarray*}
T_{3}^{\star}&=& \lambda_{n}^{2}\sum_{i =0}^{n-1}\mathbf{E}\left(K^{2}(n^{\alpha}B_{i}^{H_{1}})\int_{i}^{i+1}
 e^{i_{0}\lambda_{n}\eta_{s}^{(\lambda_{n})}}\left( \frac{1}{2} -H_{2} (s-i) ^{2H_{2}-1}\right)ds \cdot e(\lambda ,n) g_{N}\right).
\end{eqnarray*}
At this point, we will again apply the It\^o formula for $e^{i\lambda_{n}
\eta _{s}^{(\lambda_{n})}}$. It implies that
\begin{eqnarray}
T_{3}^{\star}&=& \mathbf{E}\left(\lambda_{n}^{2}\sum_{i = 0}^{n-1}K^{2}\left(n^{\alpha}B_{i}^{H_{1}}\right)\int_{i}^{i+1}
\left( \frac{1}{2} -H_{2} (s-i) ^{2H_{2}-1}\right)ds \cdot e(\lambda,n)g_{N}\right)\nonumber
\\
&&+ \mathbf{E}\left( i_{0}\lambda_{n}^{3}\sum_{i = 0}^{n-1}
K^{2}\left(n^{\alpha}B_{i}^{H_{1}}\right)\int_{i}^{i+1}\left(\frac{1}{2} -H_{2} (s-i) ^{2H_{2}-1}\right)\right.\nonumber
\\
&& \left. \int_{0}^{s} e^{i_{0}\lambda_{n}\eta_{u}^{(\lambda_{n})}}K(n^{\alpha} B^{H_{1}}_{[u]}
)dB^{H_{2}}_{u}ds \cdot e(\lambda ,n)g_{N}\right)\nonumber
\\
&&+ \mathbf{E}\left( \frac{1}{2}\lambda_{n}^{4}\sum_{i = 0}^{n-1}K^{2}\left(n^{\alpha}B_{i}^{H_{1}}\right)\int_{i}^{i+1}
\left( \frac{1}{2} -H_{2} (s-i) ^{2H_{2}-1}\right)\right. \nonumber
\\
&&\left.\int_{0}^{s} e^{i_{0}\lambda_{n} \eta_{u}^{(\lambda_{n})}}K^{2}(n^{\alpha } B^{H_{1}}_{[u]} )duds \cdot e(\lambda ,n) g_{N}\right) \nonumber
\\
&&-\mathbf{E}\left(\lambda_{n}^{4}\sum_{i = 0}^{n-1}K^{2}\left(n^{\alpha}B_{i}^{H_{1}}\right)\int_{i}^{i+1} \left( \frac{1}{2} -H_{2} (s-i) ^{2H_{2}-1}\right)\right.\nonumber
\\
&& \left.\times H_{2}(2H_{2} - 1)\int_{0}^{s} e^{i_{0}\lambda_{n} \eta_{u}^{(\lambda_{n})}}K(n^{\alpha} B_{\left[u\right]}^{H_{1}})\int_{0}^{u}K(n^{\alpha}B_{\left[v\right]}^{H_{1}})\left|u-v\right|^{2H_{2} - 2}dvduds \cdot e(\lambda, n)g_{N}du \right) \nonumber
\\
& = & b^{(1)} + b^{(2)} + b^{(3)} + b^{(4)}.\label{long}
\end{eqnarray}
The first summand $b^{(1)}$ vanishes because the integral
$$\int_{i}^{i+1}
\left( \frac{1}{2} -H_{2} (s-i) ^{2H_{2}-1}\right) ds$$ vanishes.
The second summand $b^{(2)}$ goes to zero as $n\to \infty$ using
exactly the same argument as for the convergence of $T_{1}^{\star}$.
Concerning the third summand, $b^{(3)}$, using (\ref{bound1}) and the
fact that $\vert g_{N}\vert =1$, we get
\begin{eqnarray*}
b^{(3)} & \leq & \mathbf{E}\left( \lambda_{n} ^{4}\sum_{i = 0}^{n-1}K^{2}\left(n^{\alpha}B_{i}^{H_{1}}\right)\int_{i}^{i+1}\underbrace{\left|\frac{1}{2} -H_{2} (s-i) ^{2H_{2}-1}\right|}_{\leq 1} \int_{0}^{s} K^{2}(n^{\alpha }B^{H_{1}}_{[u]})duds\right)
\\
& \leq &  \mathbf{E}\left(\frac{1}{2}\lambda_{n}^{4} \sum_{i = 0}^{n-1} K^{2}\left(n^{\alpha}B_{i}^{H_{1}}\right)\int_{i}^{i+1} \left( \sum_{j=0}^{i-1} K^{2}(n^{\alpha} B^{H_{1}}_{j}) + K^{2}(n^{\alpha} B^{H_{1}}_{i})\underbrace{(s-i)}_{\leq 1}\right)ds\right)
\\
& \leq & \mathbf{E}\left(\lambda_{n}^{4} \sum_{i = 0}^{n-1}\sum_{j=0}^{i-1}K^{2}\left(n^{\alpha}B_{i}^{H_{1}}\right)K^{2}\left(n^{\alpha}B_{j}^{H_{1}}\right)\right)
+\mathbf{E}\left(\lambda_{n}^{4} \sum_{i = 0}^{n-1}\sum_{j=0}^{i-1}K^{4}\left(n^{\alpha}B_{i}^{H_{1}}\right)
\right).
\end{eqnarray*}
The second term goes to zero because $\mathbf{E}\left(\sum_{i = 0}^{n-1}\sum_{j=0}^{i-1}K^{4}\left(n^{\alpha}B_{i}^{H_{1}}\right)\right)$ behaves as $n^{-\alpha - H_{1}+1}$ and the first term goes to zero because the non-diagonal term is dominated by the diagonal term. Analogously to the convergence of $T_{2}^{\star}$, the last summand in (\ref{long}) converges to zero. This completes the proof. \qed

\end{document}